\documentclass[a4paper,12pt]{article}
\usepackage[british]{babel}

\usepackage{rcs} \RCS $Revision: 1.7 $ \RCS $Date: 2016/04/21 22:28:10 $ \RCS $Author: hgm $ \RCS $RCSfile: Est_Parameter_Opatija.tex,v $ \RCS $Id: Est_Parameter_Opatija.tex,v 1.7 2016/04/21 22:28:10 hgm Exp $

\pdfoutput=1

\usepackage{color}

\usepackage{ucs}                
\usepackage[utf8x]{inputenc} 	

\usepackage[T1]{fontenc}     
\usepackage{lmodern}         

\usepackage[singlespacing]{setspace}

\usepackage[text={160mm,250mm},centering]{geometry}

\usepackage[sfdefault=cmbr,OMLmathsans]{isomath}  

\usepackage{authblk}
\usepackage{graphicx}
\usepackage{caption}

\usepackage{latexsym}
\usepackage{amsmath}
\usepackage{amsthm}
\usepackage{amssymb}

\usepackage{upgreek}  

\usepackage{tensor}  

\usepackage{refdef}
\usepackage{ifontdef}  

\newtheorem{thm}{Theorem}

\newcommand{\ignore}[1]{}

\newcommand{\vepsilon}{\varepsilon}

\newcommand{\vphi}{\varphi}
\newcommand{\vpi}{\varpi}

\newcommand{\vek}[1]{\mathchoice{\displaystyle\boldsymbol{#1}}
{\textstyle\boldsymbol{#1}}{\scriptstyle\boldsymbol{#1}}
{\scriptscriptstyle\boldsymbol{#1}}}

\newcommand{\ops}[1]{\mathchoice{\displaystyle\mathsf{#1}}
{\textstyle\mathsf{#1}}{\scriptstyle\mathsf{#1}}
{\scriptscriptstyle\mathsf{#1}}}

\newcommand{\tnb}[1]{\mathchoice{\displaystyle\mathboldsans{#1}}
{\textstyle\mathboldsans{#1}}{\scriptstyle\mathboldsans{#1}}
{\scriptscriptstyle\mathboldsans{#1}}}

\newcommand{\tns}[1]{\mathchoice{\displaystyle\mathsans{#1}}
{\textstyle\mathsans{#1}}{\scriptstyle\mathsans{#1}}
{\scriptscriptstyle\mathsans{#1}}}
 
\newcommand{\vcek}[1]{\vek{\check{#1}}}
\newcommand{\vbar}[1]{\vek{\bar{#1}}}
\newcommand{\vtil}[1]{\vek{\tilde{#1}}}
\newcommand{\vhat}[1]{\vek{\hat{#1}}}

\newcommand{\Tbar}[1]{\tnb{\bar{#1}}}
\newcommand{\Ttil}[1]{\tnb{\tilde{#1}}}

\newcommand{\Sbar}[1]{\tns{\bar{#1}}}
\newcommand{\Stil}[1]{\tns{\tilde{#1}}}
\newcommand{\Shat}[1]{\tns{\hat{#1}}}

\newcommand{\EXP}[1]{\mathbb{E}\left(#1\right)}

\newcommand{\diag}{\mathop{\mathrm{diag}}\nolimits}

\newcommand{\spn}{\mathop{\mathrm{span}}\nolimits}

\newcommand{\cl}{\mathop{\mathrm{cl}}\nolimits}

\newcommand{\di}{\mathrm{d}}
\newcommand{\Di}{\mathrm{D}}
\newcommand{\ii}{\mathchoice{\displaystyle\mathrm i}
{\textstyle\mathrm i}{\scriptstyle\mathrm i}
{\scriptscriptstyle\mathrm i}}

\newcommand{\citep}[1]{\cite{#1}}

\newcommand{\ip}[2]{\langle #1, #2 \rangle}

\newcommand{\nd}[1]{\| #1 \|}

\newcommand{\trpos}{\ops{T}}
\newcommand{\Id}{\mrm{I}}

\definecolor{myred}{rgb}{1, 0.2, 0.2}

\newcommand{\btheta}{\vek{\theta}}
\newcommand{\bxi}{\vek{\xi}}
\newcommand{\bbeta}{\vek{\eta}}
\newcommand{\bzeta}{\vek{\zeta}}
\newcommand{\beps}{\vek{\vepsilon}}

\newcommand{\authorhgm}{Hermann G. Matthies}
\newcommand{\authoral}{Alexander Litvinenko}
\newcommand{\authorbr}{Bojana V. Rosi\'c}
\newcommand{\authorez}{Elmar Zander}

\newcommand{\affilwire}{Institute of Scientific Computing \authorcr
                        Technische Universit\"at Braunschweig, Germany}
\newcommand{\affilkaust}{KAUST, Thuwal, Saudi Arabia}

\newcommand{\theauthor}{\authorhgm$^a$\thanks{corresponding author}\hspace{2em} 
      \authorez$^a$ \authorcr \authorbr$^a$ \hspace{2em} \authoral$^b$}

\newcommand{\thetitle}{Bayesian Parameter Estimation via\authorcr
                       Filtering and Functional Approximations}

\newcommand{\thekeywords}{inverse identification, uncertainty
  quantification, Bayesian update, parameter identification,
  conditional expectation, filters, functional and spectral approximation}

\title{\thetitle\thanks{Partly supported by the Deutsche
          Forschungsgemeinschaft (DFG) through SFB 880.}}
\author{\theauthor}

\affil{$^a$\affilwire\authorcr$^b$\affilkaust}

\date{}

\newcommand{\thebib}{./bib}

\begin{document}

\maketitle

\begin{abstract}
The \emph{inverse} problem of determining parameters in a model by
comparing some output of the model with observations is addressed.
This is a description for what hat to be done to use the 
\emph{Gauss-Markov-Kalman} filter for the \emph{Bayesian} estimation
and updating of parameters in a computational model.
This is a filter acting on random variables, and while its Monte Carlo
variant --- the \emph{Ensemble Kalman Filter} (EnKF) --- is fairly
straightforward, we  subsequently only sketch
its implementation with the help of functional representations.
\end{abstract}


\section{Introduction}  \label{S:intro}
Inverse problems in a probabilistic setting (e.g.\ \citep{jaynes03, Tarantola2004}
and references therein) are considered in here.
This situation is given in case one observes the output of some system,
and would like to infer from this the state of the system and
the values of parameters describing it, such that the output
could be caused by this combination of state and parameters.
The inverse problem is typically ill-posed, but in a probabilistic
formulation using Bayes's theorem, it becomes well-posed (e.g.\
\citep{Stuart2010}).  The unknown parameters are considered as uncertain,
and modelled as random variables (RVs).  The information available before
the measurement is called the \emph{prior}
probability distribution.  This means on one hand that the result of the
identification is a probability distribution, and not a single
value, and on the other hand the computational work may be
increased substantially, as one has to deal with RVs.  The
probabilistic setting thus can be seen as modelling our knowledge
about a certain situation --- the state and the value of the parameters --- in the
language of probability theory, and using the observation to
update our knowledge, (i.e.\ the probabilistic description) by
\emph{conditioning} on the observation.

The \emph{inverse} problem of determining or calibrating the parameters
in a computational model is addressed in the framework of Bayesian
estimation.  This is simplified to just computing the \emph{conditional expectation}.
For nonlinear models, further simplifications are needed, which
give a computationally efficient algorithm, leading via a generalisation
of the well-known \emph{Gauss-Markov} theorem to something which may be
seen as an substantial extension of the \emph{Kalman} filter.  The
resulting filter is therefore termed the \emph{Gauss-Markov-Kalman} filter (GMKF).

This document gives a short description of the connection of the
Gauss-Markov-Kalman (GMK) filter with Bayesian updating via
conditional expectation.  Subsequently it points out one of the
simplest approximations to the conditional expectation, which
results in the GMK-filter.

The key probabilistic background for this is Bayes's theorem in
the formulation of Laplace  \citep{jaynes03, Tarantola2004}.
What one wants to compute in the end are often conditional expectations
w.r.t\ the conditional distribution.  This may be achieved by
directly sampling from the conditional or posterior distribution
employing Markov-chain Monte Carlo (MCMC) methods
(see e.g.\ \citep{Hastings1970, Marzouk2007, BvrAkJsOpHgm11}).
On the other hand, it is well known
that the Bayesian update is theoretically based on the notion of
conditional expectation (CE) \citep{Bobrowski2006/087},
which may be taken as a basic theoretical notion.
It is shown that CE serves not only as a theoretical basis, but
also as a basic computational tool.  This may be seen as somewhat
related to the ``Bayes linear'' approach \citep{Goldstein2007, Kennedy01},
which has a linear approximation of CE as its basis, as will
be explained later.

In many cases, for example when tracking a dynamical system, the
updates are performed sequentially step-by-step, and for the next
step one needs not only a probability distribution in order to
perform the next step, but a random variable which may be evolved
through the state equation.  Methods on how to transform the prior RV
into the one which is conditioned on the observation will be
discussed as well \citep{HgmEzBvrAlOp15, HgmEzBvrAlOp15-p}.

The GMK-filter is so constructed that it obtains the correct posterior mean.
It is further proposed how the simple approximation of the GMK-filter,
which nevertheless is exact in some situations,
may be enhanced, so that a RV may be constructed whose distribution
approaches the posterior distribution to any desired accuracy.

This is a filter operating on random vectors, and as such needs
a \emph{stochastic discretisation} to be numerically viable.
While the numerical implementation via Monte Carlo methods ---
i.e.\ sampling or ensembles or particles --- is fairly straightforward,
here we describe the implementation via \emph{functional} approximation
or representation, where the unknown random variables are represented
as functions of \emph{known} random variables.

The plan for the rest of the paper is as follows:  in \refS{math-setup}
the mathematical set-up is described, to introduce the general setting and
mathematical background.  Finite-dimensional setting only.  But works
also in function spaces in infinite dimensions.  A synopsis of the Bayesian
approach to inverse problems is given in \refS{Bayes-synops}, stressing the
r\^ole of the conditional expectation operator (CE).  From this the filtering
approach is developed, which is described in \refS{GMKF}.  The functional
approximation is detailed in \refS{FA}, both for the filter and the forward
model.  As the Gauss-Markov-Kalman filter (GMKF) described in \refS{GMKF} is
one of the simplest but nevertheless effective approaches, some thoughts
and experiments at improved filters are given in \refS{more-acc}.
The conclusion is then given in \refS{conc}.

\section{Mathematical set-up}    \label{S:math-setup}
Assume that one has a mathematical model of the system under consideration,
symbolically written as
\begin{equation}  \label{eq:model-1}
  A(u,\vek{p}) = f,
\end{equation}
where the variable $u\in\C{U}$ represents a the \emph{state} of the
system in a vector space  $\C{U}$, the variables
$\vek{p} = [p_1,\dots,p_M] \in \C{P}=\D{R}^M \; (M\in\D{N})$
are parameters to calibrate the model,
$f\in\C{U}^*$ stands for the external influences --- the loading, action,
initial conditions, experimental set-up --- where $\C{U}^*$ is the dual
space to $\C{U}$ such that \refeq{eq:model-1} is a weak form of a state
equation, and the operator $A:\C{U}\to\C{U}^*$ describes the system under consideration.
The space $\C{U}$ may be taken as a Hilbert space for simplicity, and
later we shall assume that the model \refeq{eq:model-1} has been discretised
on some finite-dimensional subspace 
\[ \C{U}_N \subset \C{U}, \; \C{U}_N\cong\D{R}^N \; (N\in\D{N}). \]

With the help of \refeq{eq:model-1},
given an action $f\in\C{U}^*$ and a value for the
parameters $\vek{p}\in\C{P}$, we assume that it is possible to
\emph{predict} or \emph{forecast} the state $u\in\C{U}$, and from the
state it is possible to compute all other observables of the system,
see \refeq{eq:meas-1}.  In other words, the assumption is that
\refeq{eq:model-1} is well-posed, so that the state $u(\vek{p},f)$ is
a function of action $f$ and parameters $\vek{p}$.

We will tacitly assume that \refeq{eq:model-1} covers also time-evolution
problems.  To keep things notationally simple, in this case one may assume
\refeq{eq:model-1} describes the evolution over a certain time step.
The parameters $\vek{p}$ may actually include the initial conditions in case of
a time-evolution problem.

Assume also that neither the state $u\in\C{U}$ nor the parameters
$\vek{p}\in\C{P}$ are directly observable but only some function
$Y:\C{P}\times\C{U}\to\C{Y}$ of them, where the vector space
$\C{Y} \cong \D{R}^I \; (I\in\D{N})$ is assumed finite-dimensional
for the sake of simplicity.  The measurement is then
\begin{equation}  \label{eq:meas-1}
y = Y(\vek{p},u(\vek{p},f)) = Y(\vek{p}),
\end{equation}
where sometimes we shall abbreviate this simply to $y=Y(\vek{p})$ if the action
$f\in\C{U}^*$ is assumed to be given and known.

In addition there is a second system --- a more accurate one, possibly
an experiment, i.e.\ reality, something we can evaluate at possibly high cost,
but which does not need any parameters for calibration and only
serves to describe the background
\begin{equation}  \label{eq:model-2}
  A_{\circ}(u_{\circ}) = f_{\circ},
\end{equation}
where $u_\circ\in\C{U}_\circ$, again some Hilbert space not necessarily equal
to $\C{U}$ from \refeq{eq:model-1}, the right-hand-side (rhs)
$f_\circ\in\C{U}_\circ^*$ is an action, and $A_\circ:\C{U}_\circ\to\C{U}_\circ^*$.
It is assumed that $f\in\C{U}^*$ and $f_\circ\in\C{U}_\circ^*$ describe the
same situation resp.\ experiment.
Again, for the sake of simplicity, this is written in this simple
stationary form, although it may also cover evolutionary problems.

The idea is that the model in \refeq{eq:model-2} is going to be used to
calibrate --- determine the \emph{best} --- parameters $\vek{p}$ such
that the predictions of  \refeq{eq:model-1} match those of \refeq{eq:model-2}
as well as possible.
The two models --- or model and reality --- can only be compared by the
observables or measurements $y \in \C{Y}$, so we assume that there
is another function 
\begin{equation}  \label{eq:meas-2}
y_\circ = Y_{\circ}(u_{\circ}),\qquad Y_\circ:\C{U}_\circ \to \C{Y},
\end{equation}
which models the same observation in relation to \refeq{eq:model-2}.

We also assume that we observe a value $\check{y} \in \C{Y}$, which is
not directly $y_{\circ}$, but $y_{\circ}+\vepsilon$, where $\vepsilon:\Omega\to\C{Y}$
is a random variable, which in the case of \refeq{eq:model-2}
being reality models the errors of the measurement device, and in case
of \refeq{eq:model-2} being a computational
model can represent the \emph{model error} of \refeq{eq:model-2}, i.e.\ the difference
between it and \emph{reality}.  Our model for the observation of \refeq{eq:meas-2}
in terms of the quantities in \refeq{eq:model-1} is hence 
\begin{equation} \label{eq:z-err}
z = y + \vepsilon = Y(\vek{p}) + \vepsilon = Y(\vek{p},u(\vek{p},f)) + \vepsilon.
\end{equation}
This is a simple model of an additive error, which serves the purpose of
illustrating the whole procedure.
The goal of calibration is now to estimate
$\vek{p}$ such that $y$ and $y_{\circ}$ resp.\ $z$ and $\check{y}$ deviate as little
as possible.

\section{Synopsis of Bayesian estimation}   \label{S:Bayes-synops}
The idea is that the observation $z$ --- which ideally should equal $\check{y}$ ---
depends on the unknown parameters $\vek{p}$, and this should give an indication
on what $\vek{p}$ should be.  The problem in general is --- apart from the
distracting error $\vepsilon$ --- that the mapping $\vek{p} \mapsto Y(\vek{p})$
is in general not invertible, i.e.\  $z$ does not contain enough information to
uniquely determine $\vek{p}$, or there are many $\vek{p}$ which give a good
fit for $\check{y}$.  Therefore the \emph{inverse} problem of determining $\vek{p}$
from observing $\check{y}$ is termed an \emph{ill-posed} problem.

The situation is a bit comparable to Plato's allegory of the cave, where
Socrates compares the process of gaining knowledge with looking at the
shadows of the real things.  The observations $y$ resp.\ $z$ are the
``shadows'' of the ``real'' things $\vek{p}$ and $u$ resp.\ $u_\circ$, and from
observing the ``shadows'' we want to infer what ``reality'' is, in a way
turning our heads towards it.  We hence
want to ``free'' ourselves from just observing the ``shadows'' and
gain some understanding of ``reality''.

One way to deal with this difficulty is to measure the difference
between observed and predicted system output and try to find parameters such
that this difference is minimised.  Frequently it may happen that
the parameters which realise the minimum are not unique.  In case
one wants a unique parameter, a choice has to be made, usually by
demanding additionally that some norm or similar functional of the parameters
is small as well, i.e.\ some regularity is enforced. This optimisation
approach hence leads to regularisation procedures.

Here we take the view that our lack of knowledge or uncertainty
of the actual value of
the parameters can be described in a \emph{Bayesian} way through a
probabilistic model \citep{jaynes03, Tarantola2004}.  
The unknown parameter $\vek{p}$ is then modelled as a random variable
(RV)---also called the \emph{prior} model---and additional information on the
system through measurement or observation
changes the probabilistic description to the so-called \emph{posterior} model.
The second approach is thus a method to update the probabilistic description 
in such a way as to take account of the additional information, and the 
updated probabilistic description \emph{is} the parameter estimate,
including a probabilistic description of the remaining uncertainty.

It is well-known that such a Bayesian update is in fact closely related
to \emph{conditional expectation}
\citep{jaynes03, Bobrowski2006/087, Goldstein2007}, and this will be the basis
of the method presented.  For these and other probabilistic notions
see for example \citep{Papoulis1998/107} and the references therein.
As the Bayesian update may be numerically
very demanding, we show computational procedures
to accelerate this update through methods based on 
\emph{functional approximation} or \emph{spectral representation}
of stochastic problems \citep{matthies6}.  These approximations
are in the simplest case known as Wiener's so-called \emph{homogeneous}
or \emph{polynomial chaos} expansion,
which are polynomials in independent Gaussian RVs --- the ``chaos'' --- and which
can also be used numerically in a Galerkin procedure
\citep{matthies6}.

Although the Gauss-Markov theorem and its extensions \citep{Luenberger1969}
are well-known, as well as its connections to the Kalman filter
\citep{Kalman, Grewal2008} --- see also the recent Monte Carlo or
\emph{ensemble} version \citep{Evensen2009} --- the connection to Bayes's theorem is not
often appreciated, and is sketched here.   This turns out to be a linearised version
of \emph{conditional expectation} (CE).

Since the parameters of the model to be estimated are uncertain, all relevant
information may be obtained via their stochastic description.
In order to extract information from the posterior, most estimates take
the form of expectations w.r.t.\ the posterior.
These expectations --- mathematically integrals, numerically to be evaluated
by some quadrature rule --- may be computed via asymptotic,
deterministic, or sampling methods.
Here we follow our recent publications \citep{opBvrAlHgm12, BvrAkJsOpHgm11}.

To be a bit more formal, assume that the uncertain parameters are given by
\begin{equation}  \label{eq:RVp}
\vek{p}: \Omega \to  \D{R}^M \text{  as a RV on a probability space   }
  (\Omega, \F{A}, \D{P}) ,
\end{equation}
where the set of elementary events is $\Omega$, $\F{A}$ a $\sigma$-algebra of
measurable events, and $\D{P}$ a probability measure.  Additionally, also
the situation / action / loading / experiment may be uncertain, and we model
this by allowing also $f\in\C{U}^*$ in \refeq{eq:model-1} and
$f_{\circ}\in\C{U}^*_\circ$ to be random variables.
The \emph{expectation} or mean of a RV, for example $\vek{p}$,
corresponding to $\D{P}$ will be denoted by $\EXP{}$,
e.g.\ $\vbar{p}:=\EXP{\vek{p}} := \int_\Omega \vek{p}(\omega) \, \D{P}(\di \omega)$,
and the zero-mean part is denoted by $\vtil{p} = \vek{p} - \vbar{p}$.  The covariance
of $\vek{p}$ and another RV $\vek{q}$ is written as 
$\vek{C}_{pq} := \EXP{\vtil{p}\otimes\vtil{q}}$, and for short $\vek{C}_{p}$ if
$\vek{p}=\vek{q}$.

  Modelling our lack-of-knowledge
about $\vek{p}\in\C{P}$ and $u\in\C{U}$ in a Bayesian way
\citep{jaynes03, Tarantola2004, Goldstein2007} by replacing them
with  random variables (RVs), the problem becomes well-posed
\citep{Stuart2010}.  But of course one is looking now at the
problem of finding a probability distribution that best fits the data;
and one also obtains a probability distribution, not just \emph{one} pair
$\vek{p}$ and $u$.  Here we focus on the use of a linear Bayesian approach 
\citep{Goldstein2007} in the framework of ``white noise'' analysis,
but will also show some possibilities to obtain more accurate estimates
beyond the linear Bayesian approximation.

As formally $\vek{p}$ and possibly $f$ are RVs, so is the state $u(\vek{p},f)$,
reflecting the uncertainty about the state of \refeq{eq:model-1}.  From this
follows that also the prediction of the ``true'' measurement $y$ \refeq{eq:meas-1}
is a RV.
Also assume that the error $\vepsilon(\omega)$ is a $\C{Y}$-valued RV,
and in total the prediction of the observation or measurement \refeq{eq:z-err}
$z(\omega)= y(\omega) + \vepsilon(\omega)$ therefore 
becomes a RV as well; i.e.\ we have a \emph{probabilistic} description of
the prediction of the measurement.

\subsection{The theorem of Bayes and Laplace}      \label{SS:Bayes-Laplace-thm}
Bayes's original statement of the theorem which today bears his name was
only for a very special case.  The form which we know today is due to
Laplace, and it is a statement about conditional probabilities.

Bayes's theorem is commonly accepted as a consistent way to incorporate
new knowledge into a probabilistic description \citep{jaynes03, Tarantola2004}.
The elementary textbook statement of the theorem is about
conditional probabilities
\begin{equation}  \label{eq:iII}
 \D{P}(\C{I}_p|\C{M}_z) = \frac{\D{P}(\C{M}_z|\C{I}_p)}{\D{P}(\C{M}_z)}\D{P}(\C{I}_p),
 \quad \text{if }\D{P}(\C{M}_z)>0,
\end{equation}
where $\C{I}_p\subset\C{P}$ is some subset of possible $\vek{p}$'s on which we would like
to gain some information, and $\C{M}_z\subset\C{Y}$ is the information
provided by the measurement.  The term $\D{P}(\C{I}_p)$ is the so-called
\emph{prior}, it is what we know before the observation $\C{M}_z$.
The quantity $\D{P}(\C{M}_z|\C{I}_p)$ is the so-called \emph{likelihood},
the conditional probability of $\C{M}_z$ assuming that $\C{I}_p$ is given.
The term $\D{P}(\C{M}_z)$ is the so called \emph{evidence}, the probability
of observing $\C{M}_z$ in the first place, which sometimes can be expanded
with the \emph{law of total probability}, allowing to choose between
different models of explanation.  It is necessary to make the
right hand side of \refeq{eq:iII} into a real probability---summing
to unity---and hence the term $\D{P}(\C{I}_p|\C{M}_z)$, the \emph{posterior}
reflects our knowledge on $\C{I}_p$ \emph{after} observing $\C{M}_z$.

This statement \refeq{eq:iII} runs into problems if the set observations
$\C{M}_z$ has vanishing measure, $\D{P}(\C{M}_z)=0$, as is the case when we observe
\emph{continuous} random variables, and the theorem would have to be
formulated in \emph{densities}, or more precisely in 
probability density functions (pdfs).  But the statement then has the indeterminate
term $0/0$, and some form of limiting procedure is needed.  As a sign that
this is not so simple --- there are different and inequivalent forms of doing
it ---  one may just point to the so-called \emph{Borel-Kolmogorov} paradox.

There is one special case where something resembling \refeq{eq:iII} may be
achieved with pdfs, namely if $z$ and $\vek{p}$ have a 
\emph{joint} pdf $\pi_{z,p}(z,\vek{p})$.
As $z$ is essentially a function of $\vek{p}$, this is a special case
depending on conditions on the error term $\vepsilon$.  In this case
of a joint pdf Bayes's theorem \refeq{eq:iII} may be formulated as 
\begin{equation}  \label{eq:iIIa}
 \pi_{p|z}(\vek{p}|z) = \frac{\pi_{z,p}(z,\vek{p})}{Z_s(z)},
\end{equation}
where $\pi_{p|z}(\vek{p}|z)$ is the \emph{conditional} pdf, and $Z_s$
(from German \emph{Zustandssumme})
is a normalising factor such that the conditional pdf $\pi_{p|z}(\cdot|z)$
integrates to unity
\[ Z_s(z) = \int_\Omega  \pi_{z,p}(z,\vek{p}) \, \di \vek{p} . \]
The joint pdf may be split into the \emph{likelihood density} $\pi_{z|p}(z|\vek{p})$
and the \emph{prior} pdf $\pi_p(\vek{p})$
\[ \pi_{z,p}(z,\vek{p}) =  \pi_{z|p}(z|\vek{p}) \pi_p(\vek{p}) . \]
so that \refeq{eq:iIIa} has its familiar form (\citep{Tarantola2004} Ch.\ 1.5)
\begin{equation}  \label{eq:iIIa1}
 \pi_{p|z}(\vek{p}|z) = \frac{\pi_{z|p}(z|\vek{p})}{Z_s(z)}  \pi_p(\vek{p}) ,
\end{equation}
These terms are in direct correspondence with
those in \refeq{eq:iII} and carry the same names.
Once one has the conditional measure $\D{P}(\cdot|\C{M}_z)$
or even a conditional pdf $\pi_{p|z}(\cdot|z)$,
the \emph{conditional expectation} (CE) $\EXP{\cdot|z}$ may be
defined as an integral over that
conditional measure resp.\ the conditional pdf.  Thus classically,
the conditional measure or pdf implies the conditional expectation.

Please observe that
the model for the RV representing the error $\vepsilon(\omega)$ determines
the likelihood functions $\D{P}(\C{M}_z|\C{I}_q)$ resp.\ the existence and form
of the likelihood density $\pi_{z|p}(z|\vek{p})$.  
In computations, it is here that the computational
model \refeq{eq:model-1} is needed, to predict the measurement RV $z$ given
the parameters $\vek{p}$ as a RV. 

Most computational approaches determine the pdfs \citep{Hastings1970, Tarantola2004,
Marzouk2007, Stuart2010, xiu2010numerical, BvrAkJsOpHgm11, Sullivan2015},
but we will later argue that
it may be advantageous to work directly with RVs, and not with conditional
probabilities or pdfs.  To this end, the concept of conditional expectation
and its relation to Bayes's theorem is needed \citep{Bobrowski2006/087}.

\subsection{Conditional expectation}        \label{SS:CE}
To avoid the difficulties with conditional probabilities like in the
Borel-Kolmogorov paradox, \emph{Kolmogorov} himself---when he was setting up
the axioms of the mathematical theory probability---turned the relation between
conditional probability or pdf and conditional expectation around, and
defined as a first and fundamental notion \emph{conditional expectation}
\citep{Bobrowski2006/087}.

It has to be defined not with respect to measure-zero observations of a RV $z$,
but w.r.t.\ sub-$\sigma$-algebras $\F{B}\subset\F{A}$ of the underlying
$\sigma$-algebra $\F{A}$.  The $\sigma$-algebra may be loosely seen as the
collection of subsets of $\Omega$ on which we can make statements about their
probability, and for fundamental mathematical reasons in many case this is
\emph{not} the set of \emph{all} subsets of $\Omega$.  The sub-$\sigma$-algebra $\F{B}$
may be seen as the sets on which we learn something through the observation.

The simplest---although slightly restricted---way to define the conditional
expectation \citep{Bobrowski2006/087} is to just consider RVs with 
\emph{finite variance}, i.e.\ the Hilbert-space 
\[ \C{S} := L_2(\Omega,\F{A},\D{P}) := \{r:\Omega\to\D{R}\;:\;
  r \;\text{measurable w.r.t.}\;\F{A}, \EXP{|r|^2}<\infty \};\]
with the inner product given by
\begin{equation} \label{eq:ip-S}
  \forall\; r_1, r_2 \in \C{S}:\quad \ip{r_1}{r_2}_\C{S} := \EXP{r_1 \; r_2},
\end{equation}
and the usual Hilbert norm $\nd{r}_\C{S}:= \sqrt{\ip{r}{r}_\C{S}}$.
If $\F{B}\subset\F{A}$ is a sub-$\sigma$-algebra, the space
\[ \C{S}_\F{B} := L_2(\Omega,\F{B},\D{P}) := \{r:\Omega\to\D{R}\;:\;  
  r \;\text{measurable w.r.t.}\;\F{B}, \EXP{|r|^2}<\infty \} \subset \C{S}\]
is a \emph{closed} subspace, and hence has a well-defined continuous orthogonal
projection $P_\F{B}: \C{S}\to\C{S}_\F{B}$.  The \emph{conditional expectation} (CE)
of a RV $r\in\C{S}$ w.r.t.\ a sub-$\sigma$-algebra $\F{B}$ is then defined as
that orthogonal projection
\begin{equation} \label{eq:def-ce}
   \EXP{r|\F{B}} :=   P_\F{B}(r) \in \C{S}_\F{B}.
\end{equation}
It can be shown \citep{Bobrowski2006/087} to coincide with the classical notion
when that one is defined, and the \emph{unconditional} expectation $\EXP{}$
is in this view just the CE w.r.t.\ the minimal $\sigma$-algebra 
$\F{B}=\{\emptyset, \Omega\}$.
As the CE is an orthogonal projection, it minimises the squared error
\begin{equation} \label{eq:min-ce}
 \EXP{|r - \EXP{r|\F{B}}|^2} = \min\{ \EXP{|r - \hat{r}|^2}\;:\; 
   \hat{r}\in\C{S}_\F{B} \}.
\end{equation}
The CE is therefore a form of a \emph{minimum mean square error} (MMSE) estimator.
One has a form of \emph{Pythagoras's} theorem 
\[ \EXP{|r|^2} = \EXP{|\EXP{r|\F{B}}|^2} + \EXP{|r - \EXP{r|\F{B}}|^2} . \]
corresponding to the orthogonal decomposition
\[
  r = P_\F{B}(r) + (\Id_\E{S} - P_\F{B})(r) = P_\F{B}(r) +
  \left(r -  P_\F{B}(r) \right).
\]
From which --- or \refeq{eq:min-ce} --- one obtains the
\emph{variational equation} or orthogonality relation
\begin{equation} \label{eq:var-ce}
\forall \hat{r}\in\C{S}_\F{B}: \quad \EXP{\hat{r}\, (r - \EXP{r|\F{B}})}=
     \ip{\hat{r}}{r - P_\F{B}(r)}_{\C{S}} =0 .
\end{equation}

Given the CE, one may completely characterise the \emph{conditional} probability,
e.g.\ for $\C{A} \subset \Omega, \C{A} \in \F{B}$ by 
\[ \D{P}(\C{A} | \F{B}) := \EXP{\chi_{\C{A}}|\F{B}} , \]
where $\chi_{\C{A}}$ is the RV which is unity iff $\omega\in\C{A}$ and vanishes otherwise
--- the \emph{usual} characteristic function, sometimes also termed an indicator
function.  Thus if we know $\D{P}(\C{A} | \F{B})$ for each $\C{A} \in \F{B}$, we know
the conditional probability.  Hence having the CE $\EXP{\cdot|\F{B}}$ allows one
to know everything about the conditional probability.
If the prior probability was the distribution of some
RV $r$, we know that is is completely characterised by the \emph{prior} characteristic
function --- in the sense of probability theory --- $\vphi_r(s) := \EXP{\exp(\ii r s)}$.
To get the \emph{conditional} characteristic function $\vphi_{r|\F{B}}(s) =
\EXP{\exp(\ii r s)|\F{B}}$, all one has to do is use the CE instead of the unconditional
expectation.  This then completely characterises the conditional distribution.

In our case of an observation of a RV $z$,
the sub-$\sigma$-algebra $\F{B}$ will be the one generated by
the \emph{observation} $z$, i.e.\ $\F{B}=\sigma(z)$, these are those
subsets of $\Omega$ on which we may obtain \emph{information} from the
observation.  According to the \emph{Doob-Dynkin} lemma the subspace
$\C{S}_{\sigma(z)}$ is given by
\begin{equation} \label{eq:DDL}
  \C{S}_{\sigma(z)} :=  \{r \in \C{S} \;:\; r(\omega) = \phi(z(\omega)),  
  \phi \;\text{measurable} \} \subset \C{S} ,
\end{equation}
i.e.\ functions of the observation.  This means intuitively that anything we
learn from an observation is a function of the observation, and the subspace
$\C{S}_{\sigma(z)} \subset \C{S} $ is where the information from the measurement
is lying.

As according to \refeq{eq:def-ce} $\EXP{r|{\sigma(z)}}=P_{\sigma(z)}(r)
\in \C{S}_{\sigma(z)}$, it is clear from \refeq{eq:DDL} that there is a measurable
function $\vpi_r$ on $\C{Y}$ such that
\begin{equation} \label{eq:DDL-csx}
  \EXP{r|{\sigma(z)}}=P_{\sigma(z)}(r) = \vpi_r(z).
\end{equation}
Observe that the CE $\EXP{r|\sigma(z)}$ and conditional probability
$\D{P}(\C{A}|\sigma(z))$---which we will abbreviate to
$\EXP{r|z}$, and similarly $\D{P}(\C{A} | \sigma(z))=\D{P}(\C{A}|z)$---is a RV,
as $z(\omega)$ is a RV.  Once an observation has been made,
i.e.\ we observe for the RV $z$ the fixed value $\check{y}\in\C{Y}$,
then $\EXP{r|\check{y}} \in \D{R}$ is just a number --- the \emph{posterior expectation},
and $\D{P}(\C{A}|\check{y})=\EXP{\chi_{\C{A}}|\check{y}}$ --- for almost all
$\check{y}\in\C{Y}$ --- is the \emph{posterior probability}.
Often these are also termed conditional expectation and conditional probability,
which may lead to confusion.  We therefore prefer the attribute \emph{posterior}.

In relation to Bayes's theorem, one may conclude that if it is possible to
compute the CE w.r.t.\ an observation $z$ or rather the posterior expectation,
then the conditional and especially the posterior probabilities after the
observation $z=\check{y}$ may as well e computed, regardless of the case whether joint
pdfs exist or not.  We take this as the starting point to Bayesian estimation.

The conditional expectation has been formulated for scalar RVs, but it is clear
that the notion carries through to vector-valued RVs in a straightforward manner,
formally by seeing a --- let us say --- $\C{V}$-valued RV as an element of the
tensor Hilbert space $\E{V}=\C{V}\otimes\C{S}$.   The CE on $\E{Y}$ is then formally
given by $\D{E}_{\E{V}}(\cdot|\F{B}):=\Id_{\C{V}}\otimes\EXP{\cdot|\F{B}}=
\Id_{\C{V}}\otimes P_{\F{B}}$, where
$\Id_{\C{V}}$ is the identity operator on $\C{V}$.  It is an orthogonal projection
in $\E{V}$, for simplicity also denoted by $P_{\F{B}}$ (cf.\ \refeq{eq:ip-Ht}).

\section{The Gauss-Markov-Kalman filter (GMKF)}       \label{S:GMKF}
It turned out that practical computations in the context of Bayesian estimation
can be extremely demanding, see \citep{McGrayne2011} for an account of the
history of Bayesian theory, and the break-throughs required in computational
procedures to make Bayesian estimation possible at all for practical purposes.
This involves both the Monte Carlo (MC) method and the Markov chain Monte Carlo
(MCMC) sampling procedure.

To arrive at computationally feasible procedures for computationally
demanding models \refeq{eq:model-1}, where MCMC methods are not feasible,
approximations are necessary.  This means in some way not using all information
but having a simpler computation.  Incidentally, this connects with the
Gauss-Markov theorem \citep{Luenberger1969} and the Kalman filter (KF)
\citep{Kalman, Grewal2008}.  These were initially stated and developed without
any reference to Bayes's theorem.  The Monte Carlo (MC) computational
implementation of this is the \emph{ensemble} KF (EnKF) \citep{Evensen2009}.
We will in contrast use a white noise or polynomial chaos approximation
\citep{opBvrAlHgm12, BvrAkJsOpHgm11}.  But the initial ideas leading to
the abstract Gauss-Markov-Kalman filter (GMKF) are independent of any
computational implementation and are presented first.  It is in an abstract
way just \emph{orthogonal projection}.

\subsection{Orthogonal decomposition}      \label{SS:orth-dec}
Assuming that the Hilbert space $\C{V}$ has inner product $\ip{\cdot}{\cdot}_\C{V}$, one
defines the Hilbert tensor inner product for elementary tensors $v\otimes r \in 
\E{V}=\C{V}\otimes\C{S}$ by
\begin{equation} \label{eq:ip-Ht}
  \forall\; v_1\otimes r_1, v_2\otimes r_2 \in \E{V}=\C{V}\otimes\C{S}:\quad 
     \ip{v_1\otimes r_1}{v_2\otimes r_2}_\E{V} := \ip{v_1}{v_2}_\C{V} \cdot 
     \ip{r_1}{r_2}_\C{S},
\end{equation}
and extends this to all of $\E{V}=\C{V}\otimes\C{S}$ by linearity.  This then
defines the Hilbert norm $\nd{v\otimes r}_\E{V}:= \sqrt{\ip{v\otimes r}{v\otimes r}_\E{V}}
= \nd{v}_\C{V} \cdot \nd{r}_\C{S}$ 
on $\E{V}$ in the usual way.  Given two RVs $\tns{v}_1, \tns{v}_2 \in 
\E{V}=\C{V}\otimes\C{S}$, they are called \citep{Bosq2000} \emph{weakly orthogonal}
iff $\ip{\tns{v}_1}{\tns{v}_2}_\E{V}=0$;

Considering now a subspace $\E{V}_\F{B} := \C{V}\otimes\C{S}_\F{B}$ with orthogonal
projector $P_\F{B}$, a RV $\tns{v}\in \E{V}$ may be decomposed into its orthogonal
components w.r.t.\ $\E{V}_\F{B}$ by
\begin{equation} \label{eq:orth-comp}
  \tns{v} = P_\F{B}(\tns{v}) + (\Id_\E{V} - P_\F{B})(\tns{v}) = P_\F{B}(\tns{v}) +
  \left(\tns{v} -  P_\F{B}(\tns{v}) \right),
\end{equation}
where $(\Id_\E{V} - P_\F{B})(\tns{v})\in \E{V}_\F{B}^\perp$, the orthogonal complement
of $\E{V}_\F{B}$.  Obviously $P_\F{B}(\tns{v})$ is the best estimator for $\tns{v}$
--- measured in the error norm squared $\nd{\tns{v}-P_\F{B}(\tns{v})}_\E{V}^2$ --- from
the subspace $\E{V}_\F{B}$.  Obviously, analogous to \refeq{eq:var-ce}, one has
\begin{equation} \label{eq:var-ce-vek}
  \forall \Shat{v}\in\E{V}_\F{B}:\qquad \ip{\Shat{v}}{\tns{v} -  P_\F{B}(\tns{v})}_{\E{V}}=0 .
\end{equation}

Further the notion of \emph{correlation} and \emph{covariance}
operator is needed:  Given a RV $\tns{v}_1\in\E{V}_1=\C{V}_1\otimes\C{S}$,
and a RV $\tns{v}_2\in\E{V}_2=\C{V}_2\otimes\C{S}$, their
correlation operator $\hat{C}_{v_1 v_2} \in \E{L}(\C{V}_2,\C{V}_1)$ --- a linear
operator from $\C{V}_2$ to $\C{V}_1$ --- is given for $w\in\C{V}_2$ by
$\hat{C}_{v_1 v_2}(w) = \D{E}_{\C{V}_1}( \ip{\tns{v}_2}{w}_{\C{V}_2} \cdot \tns{v}_1)
\in \C{V}_1$.   
The \emph{covariance} operator $C_{v_1 v_2} \in \E{L}(\C{V}_2,\C{V}_1)$
is defined by looking at the \emph{zero-mean} versions of the RVs,
i.e.\  $\Stil{v}:= \tns{v}-\EXP{\tns{v}}=\tns{v}-\Sbar{v}$.  The \emph{covariance}
operator $C_{v_1 v_2}$ is then the correlation operator of the zero-mean RVs
$\Stil{v}_1$ and $\Stil{v}_2$.  In case $\tns{v}=\tns{v}_1=\tns{v}_2$,
we for brevity just write $\hat{C}_v$ and $C_v$.

Two vector-valued RVs $\tns{v}_1, \tns{v}_2 \in \E{V}=\C{V}\otimes\C{S}$ are
called \emph{orthogonal} iff $\ip{\tns{v}_1}{(L\otimes \Id_\C{S}) \tns{v}_2}_\E{V}=0$ for all
$L \in \E{L}(\C{V},\C{V})$.  Obviously orthogonality implies weak orthogonality,
but not the other way around.  Two zero-mean RVs are called \emph{uncorrelated} iff
they are orthogonal, and in this case their covariance operator $C_{v_1 v_2}=
\EXP{\tns{v}_1\otimes\tns{v}_2}$ vanishes.

A subspace $\E{V}_s \subset \E{V} = \C{V}\otimes\C{S}$
is called $\E{L}$-\emph{closed} \citep{Bosq2000}, iff 
\[ \E{V}_s = \{ (L\otimes \Id_\C{S}) \tns{v}\;:\;  L\in\E{L}(\C{V},\C{V}),
   \;\tns{v} \in \E{V}_s \} , \]
i.e.\ $\E{V}_s$ is \emph{invariant} under all linear maps in
$L\in\E{L}(\C{V},\C{V})$.  A subspace
$\E{V}_s \subset \E{V}$ is called \emph{zero-mean} iff all RVs in $\E{V}_s$
are zero mean, and the notions of weak orthogonality and orthogonality can be
extended to subspaces in a natural fashion.

Looking at subspaces of the form $\E{V}_\F{B}=\C{V}\otimes\C{S}_\F{B}$ considered
previously, it is clear that they are $\E{L}$-closed, and hence the decompositions in
\refeq{eq:orth-comp}, \refeq{eq:orth-best}, and  \refeq{eq:orth-CE} are not
just weakly orthogonal but
in fact orthogonal in the sense just explained, i.e.\ under $\E{L}$-invariance.
This means that in addition to \refeq{eq:var-ce-vek} one has the stronger condition
\begin{equation} \label{eq:var-ce-corr}
  \forall \Shat{v}\in\E{V}_\F{B}:\qquad \hat{C}_{\hat{v},\tns{v} -  P_\F{B}(\tns{v})}=
  \EXP{\Shat{v}\otimes(\tns{v} -  P_\F{B}(\tns{v}))}=0 .
\end{equation}

\subsection{Building the filter}          \label{SS:build-filt}
Reverting to the problem of estimation after a measurement $z$, we see that
all operations are performed as projections in vector spaces.  It is possible
that the parameters $\vek{p}$ are not \emph{free} in a vector space, i.e.\ some
components may be required to be positive, or lie between some finite bounds.

This is detrimental for the estimation process, and we transform the
parameters with an invertible transformation $X:\C{X}\to\C{P}$ by
$\vek{p}=X(\vek{x})$ --- $(\vek{x}\in\C{X}=\D{R}^M)$ ---
onto \emph{new} parameters which will be estimated, and which have \emph{no}
constraints.

For simplicity one may assume that for $m=1,\dots,M$ one has $p_m = X_m(x_m)$,
for $\vek{x}\in\D{R}^M = \C{X}$.  In case $p_m$ is constrained to an
one-sided semi-infinite interval,
e.g.\ positivity, the \emph{logarithm} and its inverse the \emph{exponential}
function can be used for the transform
$X_m$, after a possible shift and scaling.  Similarly, if $p_m$ is constrained
to a finite interval, after a possible shift and scaling the \emph{probit}, \emph{logit},
or \emph{arctan} and their well-known inverse functions can be used for the
transform.

Assuming that this has been carried out, we consider now the problem of
estimating $\vek{x}$.  To be concrete, assume also that there are $I\in\D{N}$
measurements, i.e.\ the total measurement $\vek{y}\in \C{Y}$ in
\refeq{eq:meas-1} lies in $\C{Y} = \D{R}^I$, an $I$-dimensional space.
For \refeq{eq:meas-1} we shall now just write for simplicity by abuse of notation
interchangably
\begin{equation}   \label{eq:meas-v}
   \vek{y} = Y(\vek{x}) = Y(X(\vek{x})) = Y(\vek{p}).
\end{equation}

\subsubsection{The conditional expectation mean filter}          \label{SSS:CE-filt}
Now reverting to
the problem of estimating $\vek{x}_a$ from a forecast $\vek{x}_f$ and an
observation $\vcek{y}$ for the RV $\vek{z} = Y(\vek{x}) + \vek{\vepsilon}$, we consider
the subspace $\E{X}_{\sigma(z)} = \C{X}\otimes \C{S}_{\sigma(z)} \subset
\E{X} = \C{X}\otimes \C{S}$, and for $\tnb{v} = \Psi(\tnb{x})$ in \refeq{eq:orth-comp}
which can be any measurable function $\Psi$ of $\vek{x}(\omega)$ --- which in
the tensor product $\E{X}$ is denoted by $\tnb{x}\in\E{X}$ --- we take the
identity $\Psi(\tnb{x}) = \tnb{x}$ in \refeq{eq:orth-comp}.   The orthogonal
decomposition \refeq{eq:orth-comp} is for this case
\begin{equation} \label{eq:orth-best}
  \tnb{x} = P_{\sigma(z)}(\tnb{x}) + (\Id_\E{X} - P_{\sigma(z)})(\tnb{x}) =
  \D{E}_\E{X}(\tnb{x}|z) + \left(\tnb{x} - \D{E}_\E{X}(\tnb{x}|z) \right),
\end{equation}
with $P_{\sigma(z)}(\Id_\E{X} - P_{\sigma(z)})(\tnb{x}) = 0$.
Now let $\tnb{x}_f$ be the forecast --- i.e.\ representing our knowledge
before the observation $\vek{z}=\vcek{y}$ --- to which \refeq{eq:orth-best}
applies as well.   An observation $\vcek{y}$  will tell us something about the
first component in \refeq{eq:orth-best}.  Hence one defines
the \emph{filtered}, \emph{analysed}, or \emph{assimilated} RV
$\tnb{x}_a$ \emph{after} the observation $\vcek{y}$
\begin{equation} \label{eq:orth-CE}
  \tnb{x}_a = 
  \D{E}_\E{X}(\tnb{x}_f|\vcek{y}) + \left(\tnb{x}_f - \D{E}_\E{X}(\tnb{x}_f|z) \right)=
  \tnb{x}_f + \left( \D{E}_\E{X}(\tnb{x}_f|\vcek{y}) - \D{E}_\E{X}(\tnb{x}_f|z) \right)
  = \tnb{x}_f + \tnb{x}_i ,
\end{equation}
where $\tnb{x}_i= \D{E}_\E{X}(\tnb{x}_f|\vcek{y}) - \D{E}_\E{X}(\tnb{x}_f|z)$ is
called the \emph{innovation}.  This means the first term  in \refeq{eq:orth-best}
has been changed to the posterior CE $\D{E}_\E{X}(\tnb{x}_f|\vcek{y})$, and the
rest, the second term in \refeq{eq:orth-best}, has been left unchanged.  The
\refeq{eq:orth-CE} is called the \emph{conditional expectation mean} filter (CEMF).
The following is an easy consequence of the previous development:

\begin{thm}  \label{T:CEF}
With the shorthand
$\Tbar{x}^{\vcek{y}} = \D{E}_\E{X}(\tnb{x}|\vcek{y})$ and $\Ttil{x}^{\vcek{y}} =
\tnb{x} - \Tbar{x}^{\vcek{y}}$, one has
\begin{equation} \label{eq:CE-rels}
  \Tbar{x}^{\vcek{y}}_a = \Tbar{x}^{\vcek{y}}_f = \D{E}_\E{X}(\tnb{x}_f|\vcek{y}),
  \qquad \text{and} \qquad  \Ttil{x}^{\vcek{y}}_a =
   ( \tnb{x}_f - \D{E}_\E{X}(\tnb{x}_f|z) )
\end{equation}
for the posterior mean and the posterior zero-mean part.
The RV $\Ttil{x}^{\vcek{y}}$ in \refeq{eq:CE-rels} is \emph{uncorrelated} and hence
$\E{L}$-orthogonal to all RVs in $\E{X}_{\sigma(z)}$, and from \refeq{eq:orth-best}
it follows that $\D{E}_\E{X}(\tnb{x}_a|\vcek{y})
=\D{E}_\E{X}(\tnb{x}_f|\vcek{y})$.  In other
words, $\tnb{x}_a$ is \emph{unbiased}, 
$\D{E}_\E{X}(\tnb{x}_f|\vcek{y})$ is optimal,
i.e.\ the \emph{best unbiased estimator}, and $\Ttil{x}^{\vcek{y}}$ is the orthogonal error.
The posterior covariance of
$\tnb{x}_a$ is thus
\begin{equation} \label{eq:CE-cov}
   \vek{C}_{x_a} = \EXP{\vtil{x}^{\vcek{y}}_a \otimes \vtil{x}^{\vcek{y}}_a | \vcek{y}} .
\end{equation}

The RV $\tnb{x}_a$ has thus the same posterior expected value as
the posterior Bayesian distribution after the observation $\vcek{y}$, and its 
posterior covariance is given by \refeq{eq:CE-cov}, which is general \emph{not} the
correct covariance of the posterior Bayesian distribution.

It is well known \citep{Tarantola2004, Goldstein2007, Bosq2000}
that in case the prior or forecast RV $\tnb{x}_f$ is Gaussian,
the observation $Y(\vek{x}_f)$ is linear in $\tnb{x}_f$, and the error $\vepsilon$
also Gaussian, that the distribution of $\tnb{x}_a$ is the \emph{exact}
Bayesian posterior, and not just its expected value equal to the Bayesian
mean.
\end{thm}

From the fact that $\D{E}_\E{X}(\tnb{x}_f|z) = P_{\sigma(z)}(\tnb{x}_f) \in
\E{X}_{\sigma(z)}$, and knowing from the Doob-Dynkin lemma \refeq{eq:DDL} that
\begin{equation} \label{eq:X-full}
\E{X}_{\sigma(z)} = \{ \tnb{w} \in \E{X} \; :\; \tnb{w} = \phi(\vek{z}(\omega)), \;
 \phi:\C{Y}\to\C{X} \quad \text{measurable} \}, 
\end{equation}
it is clear from \refeq{eq:DDL-csx} that there is a
measurable map $\vpi_\Psi:\C{Y}\to\C{X}$ such that
\begin{equation} \label{eq:opt-map}
  \D{E}_\E{X}(\tnb{x}_f|\vek{z}) = P_{\sigma(z)}(\tnb{x}_f) = \vpi_\Psi(\vek{z}).
\end{equation}
The optimal map $\vpi_\Psi$ obviously depends on the function $\Psi(\tnb{x})$ for
which it is determined.  As we are here interested in $\Psi(\tnb{x})=\tnb{x}$,
we shall denote the optimal map in this case by $\vpi_x$.
From \refeq{eq:min-ce} one may show that $\vpi_x$ is defined by
\begin{equation} \label{eq:min-CE}
  \nd{\tnb{x}_f - \vpi_x(\vek{z})}^2_{\E{X}} = 
    \min_{\phi} \nd{\tnb{x}_f - \phi(\vek{z})}^2_{\E{X}}
  = \min_{\tnb{w}\in\E{X}_{\sigma(z)}} \nd{\tnb{x}_f - \tnb{w}}^2_{\E{X}},
\end{equation}
where $\phi$ ranges over all measurable maps $\phi:\C{Y}\to\C{X}$.  Observe that
although the minimising point $\vpi_x(\vek{z})$ is unique, the map $\vpi_x:\C{Y}\to\C{X}$
may not necessarily be so.

As $\E{X}_{\sigma(z)}$ is $\E{L}$-closed, it is characterised as in
\refeq{eq:var-ce-corr} by orthogonality in the $\E{L}$-invariant sense
\begin{equation} \label{eq:var-CE}
  \forall \tnb{w}\in\E{X}_{\sigma(z)}:\quad \hat{C}_{w, (x_f - \vpi_x(z))} 
  =\EXP{\tnb{w}\otimes(\tnb{x}_f - \vpi_x(\vek{z}))}= 0,
\end{equation}
i.e.\ the RV $(\tnb{x}_f - \vpi_x(\tnb{z}))$ is orthogonal in the $\E{L}$-invariant sense
to all RVs $\tnb{w}\in\E{X}_{\sigma(z)}$.  This \refeq{eq:var-CE} is the relation which
is used to determine $\vpi_x$.

The  assimilated RV $\tnb{x}_a$ after the observation $\vcek{y}$ in \refeq{eq:orth-CE}
is thus given by the CEM-filter equation
\begin{equation} \label{eq:orth-CE-map}
  \tnb{x}_a =   \tnb{x}_f + ( \vpi_x(\vcek{y}) - \vpi_x(\tnb{z}) )  = \tnb{x}_f + \tnb{x}_i .
\end{equation}
The terms in \refeq{eq:CE-rels} are hence given by
\begin{equation} \label{eq:CE-rels-map}
  \Tbar{x}^{\vcek{y}}_a = \vpi_x(\vcek{y}), \qquad \text{and} \qquad 
  \Ttil{x}^{\vcek{y}}_a = \tnb{x}_f - \vpi_x(\tnb{z}) .
\end{equation}

\refeq{eq:orth-CE-map} is the \emph{best unbiased} filter, with $\vpi(\check{y})$
a MMSE estimate.  Although the CE
$\D{E}_\E{X}(\tnb{x}_f|z) = P_{\sigma(z)}(\tnb{x}_f)$ is an orthogonal projection,
as the measurement operator $Y(\vek{x})$ in $\vek{z}=Y(\vek{x})+\vek{\vepsilon}$ is not
necessarily linear in $\vek{x}$,  neither is the optimal map $\vpi_x(\tnb{z})$.

\subsubsection{The linear filter}      \label{SSS:lin-filt}
The minimisation in \refeq{eq:min-CE} over all measurable maps is still a
formidable task, and typically only feasible in an approximate way.  Thus
we replace  $\E{X}_{\sigma(z)}$ by a smaller subspace; and we choose
in some way the simplest possible one
\begin{equation} \label{eq:X-1}
\E{X}_1 = \{ \tnb{w} \;:\;\tnb{w} = \phi(\tnb{z})= \vek{L} (\vek{z}(\omega)) + \vek{b}, \;
 \vek{L} \in \E{L}(\C{Y},\C{X}),\; \vek{b} \in \C{X}  \} \subset \E{X}_{\sigma(z)}
 \subset \E{X} ,
\end{equation}
where the $\phi$ are just \emph{affine} maps; they are certainly measurable.  Note that
$\E{X}_1$ is also an $\E{L}$-invariant subspace of $\E{X}_{\sigma(z)}\subset\E{X}$.
Note that also other, possibly larger, $\E{L}$-invariant subspaces of $\E{X}_{\sigma(z)}$
can be used, but this seems to be smallest useful one.
Now the minimisation \refeq{eq:min-CE} may be replaced by 
\begin{equation} \label{eq:min-CE-lin}
  \nd{\tnb{x}_f - (\vek{K}(\vek{z})+\vek{a})}^2_{\E{X}} = 
  \min_{\vek{L}, \vek{b}} \nd{\tnb{x}_f - (\vek{L}(\vek{z})+\vek{b})}^2_{\E{X}} ,
\end{equation}
and the optimal affine map is defined by $\vek{K}\in \E{L}(\C{Y},\C{X})$ and
$\vek{a} \in \C{X}$.

Using this instead of $\vpi_x$ in \refeq{eq:orth-CE-map}, one disregards some
information as $\E{X}_1 \subset \E{X}_{\sigma(z)}$ is a true subspace --- 
observe that the subspace represents the information we may learn from
the measurement --- but the computation is easier, and one arrives at
\begin{equation} \label{eq:orth-CE-lin}
  \tnb{x}_a =  \tnb{x}_f + ( \vek{K}(\vcek{y}) - \vek{K}(\vek{z}) )= 
  \tnb{x}_f + \vek{K}(\vcek{y} - \vek{z}(\omega)) = 
  \tnb{x}_f + \vek{K}(\vcek{y} - (Y(\vek{x}_f(\omega))+\vek{\vepsilon}(\omega))).
\end{equation}
This is the \emph{best linear} filter, with the linear MMSE $\vek{K}(\vcek{y})$.
One may note that the constant term
$\vek{a}$ in \refeq{eq:min-CE-lin} drops out in the filter equation.

\subsection{The Gauss-Markov theorem and the Kalman filter}  \label{SS:GMthm-KF}
The optimisation described in \refeq{eq:min-CE-lin} is a familiar one, it is
easily solved, and the solution is given by an extension of the \emph{Gauss-Markov}
theorem \citep{Luenberger1969}.  The same idea of a linear MMSE is behind the
\emph{Kalman} filter 
\citep{Kalman, Grewal2008, Goldstein2007, Papoulis1998/107, Evensen2009}.
In our context it reads
\begin{thm}  \label{T:GM}
The solution to \refeq{eq:min-CE-lin}, minimising
\[ \nd{\tnb{x}_f - (\vek{K}(\tnb{z})+\vek{a})}^2_{\E{X}} =  
  \min_{\vek{L}, \vek{b}} \nd{\tnb{x}_f - (\vek{L}(\tnb{z})+\vek{b})}^2_{\E{X}} \]
is obtained via the analog of \refeq{eq:var-CE} and
is given by $\vek{K} := \vek{C}_{x_f z} \vek{C}_{z}^{-1}$ and
$\vek{a} := \Tbar{x}_f - \vek{K}(\vbar{z})$, where $\vek{C}_{x_f z}$
is the covariance of $\tnb{x}_f$ and $\tnb{z}$, and $\vek{C}_{z}$ is the
auto-covariance of $\tnb{z}$.    In case
$\vek{C}_{z}$ is \emph{singular}, the \emph{pseudo-inverse} can be taken
instead of the inverse.
\end{thm}
The operator $\vek{K}$ is also called the \emph{K\'alm\'an} gain, and has
the familiar form known from least squares projections.  It is interesting
to note that initially the connection between MMSE and Bayesian estimation
was not seen \citep{McGrayne2011}, although it is one of the simplest
approximations to the Bayesian estimate.

The resulting filter --- with the understanding that $\vek{C}^{-1}_{z}$
is the pseudo-inverse in case of singularity ---
\begin{equation} \label{eq:GMKF}
  \tnb{x}_a =  \tnb{x}_f + \vek{C}_{x_f z} \vek{C}^{-1}_{z}(\vcek{y} - \vek{z}(\omega)) 
    = \tnb{x}_f + \vek{K}(\vcek{y} - \tnb{z}),
\end{equation}
is therefore called the \tbf{Gauss-Markov-Kalman} filter (GMKF).  Observe that in
case $\tnb{x}_f$ resp.\ $\vek{y}=Y(\tnb{x}_f)$ and $\vek{\vepsilon}$ are \emph{independent} RVs
--- as can often be assumed --- then simply $\vek{C}_{z} = \vek{C}_{y} + \vek{C}_{\vepsilon}$.

The Kalman filter has \refeq{eq:GMKF} for the means, which is obtained by
taking the expected value on both sides of \refeq{eq:GMKF}, i.e.\ due to
linearity of the expectation of each term individually:
\[ \vbar{x}_a = \vbar{x}_f + \vek{K}(\vcek{y} - \vbar{z}) . \]
It easy to compute that  \citep{Luenberger1969}
\begin{thm}   \label{T:GM-var}
The posterior covariance operator $\vek{C}_{x_a}=\EXP{\vtil{x}^{\vcek{y}}_a \otimes
\vtil{x}^{\vcek{y}}_a}$ in \refeq{eq:CE-cov} of $\tnb{x}_a$ from \refeq{eq:GMKF} is given by
\begin{equation}  \label{eq:GMKF-var}
 \vek{C}_{x_a} = \vek{C}_{x_f} - \vek{K}\vek{C}_{x_f z}^{\trpos}
    = \vek{C}_{x_f}-\vek{C}_{x_f z}\vek{C}_{z}^{-1}\vek{C}_{x_f z}^{\trpos} , 
\end{equation}
which is K\'alm\'an's formula for the covariance.
\end{thm}
This shows that \refeq{eq:GMKF} is a true extension of the classical Kalman filter (KF).
It also shows that $\vek{C}_{x_a} \leq \vek{C}_{x_f}$ in the usual ordering of 
symmetric positive definite (spd) matrices, as the spd-term
$\vek{C}_{x_f z}\vek{C}_{z}^{-1}\vek{C}_{x_f z}^{\trpos}$ is subtracted from $\vek{C}_{x_f}$.

Rewriting \refeq{eq:GMKF} explicitly in less symbolic notation
\begin{equation} \label{eq:GMKF-RV}
  \vek{x}_a(\omega) =  \vek{x}_f(\omega) + 
      \vek{C}_{x_f z} \vek{C}^{-1}_{z}(\vcek{y} - \vek{z}(\omega)) =
      \vek{x}_f(\omega) + \vek{K}(\vcek{y} - \vek{z}(\omega)) ,
\end{equation}
one may see that it is a relation between RVs, and hence some further \emph{stochastic}
discretisation is needed for it to be numerically implementable.

\section{Functional approximation}          \label{S:FA}
Our starting point is \refeq{eq:GMKF-RV}.   As it is a relation between RVs, it
certainly also holds for \emph{samples} of the RVs.  Thus it is possible to
take an \emph{ensemble} of sampling points $\omega_1,\dots,\omega_S$ and
require
\begin{equation} \label{eq:EnKF}  \forall s = 1,\dots,S:\quad
  \vek{x}_a(\omega_s) =  \vek{x}_f(\omega_s) + 
      \vek{C}_{x_f z} \vek{C}^{-1}_{z}(\check{y} - \vek{z}(\omega_s)) ,
\end{equation}
and this is the basis of the \emph{ensemble} KF, the EnKF \citep{Evensen2009};
the points $\vek{x}_f(\omega_s)$ and $\vek{x}_a(\omega_s)$ are sometimes
also denoted as \emph{particles}, and \refeq{eq:EnKF} is a simple version of
a \emph{particle filter}.

Some of the main work for the EnKF consists in obtaining
good estimates of $\vek{C}_{x_f z}$ and $\vek{C}_{z}$, as they have to be computed
from the samples.  Further approximations are possible, for example such as
\emph{assuming} a particular form for $\vek{C}_{x_f z}$ and $\vek{C}_{z}$.
This is the basis for methods like \emph{kriging} and \emph{3DVAR} resp.\ \emph{4DVAR},
where one works with an approximate Kalman gain $\vtil{K} \approx \vek{K}$.

To actually compute \refeq{eq:EnKF}, one needs to evaluate the term
$\vek{z}(\omega_s) = Y(\vek{p}(\omega_s),u(\omega_s)) + \vek{\vepsilon}(\omega_s)$
from \refeq{eq:meas-1}.  Here one the state of the system $u(\omega_s)$ being
observed and identified appears.  As alluded to in \refS{math-setup}, a 
finite dimensional approximation or \refeq{eq:model-1} is necessary for a numerical
evaluation.   This may be achieved by different means, e.g.\ finite elements, finite
volumes, finite differences, etc., e.g.\ \citep{strangFix88}.
We shall assume only that there is a finite
dimensional computational model on $\C{U}_N \subset \C{U}$ with $\C{U}_N \cong \D{R}^N$:
\begin{equation}  \label{eq:model-N-p}
  \vek{A}(\vek{u},\vek{x}) = \vek{f},
\end{equation}
an equation in $\C{U}_N^*$.  Numerical methods to solve \refeq{eq:model-N} for
$\vek{u}\in\C{U}_N$ given $\vek{f}\in\C{U}_N^*$ and $\vek{x}\in\C{X}$ typically
drive the residuum to naught:
\begin{equation}  \label{eq:resid-N-p}
  \vek{r}(\vek{u}) = \vek{f} - \vek{A}(\vek{u},\vek{x}) = 0.
\end{equation}
This now is a computational model which --- given a sample $\omega_s$ --- 
may be used to compute $\vek{u}(\omega_s)$ from $\vek{r}(\vek{u}(\omega_s))=0$,
to be used in the evaluation of
$\vek{z}(\omega_s) = Y(\vek{x}(\omega_s),\vek{u}(\omega_s)) + \vek{\vepsilon}(\omega_s)$
above.  Now all terms in \refeq{eq:EnKF} can be evaluated.

\subsection{Basics of functional approximation}      \label{SS:FA-basic}
Here we want to pursue a different tack, and want to discretise RVs not through
their samples $\omega_s$, but by \emph{functional approximations} 
\citep{matthies6, xiu2010numerical, knio2010spectral}.  This means that all RVs,
say $\vek{u}(\omega)$, are described as functions of \emph{known} RVs
$\{\xi_1(\omega),\dots,\xi_n(\omega),\dots\}$.  Often, when for example
stochastic processes or random fields are involved, one has to deal here
with \emph{infinitely} many RVs, which for an actual computation have to be
truncated to a finte vector $\bxi(\omega)=[\xi_1(\omega),\dots,\xi_L(\omega)]\in
\Xi \cong \D{R}^L$ of significant RVs.
We shall assume that these have been chosen such as to be
independent, and often even normalised Gaussian and independent.  We shall
assume that these are used to describe both the uncertainty in the parameters
$\vek{x}=[x_1,\dots,x_M]^{\trpos}$ as well as in the rhs $\vek{f}=[x_1,\dots,x_N]^{\trpos}$,
i.e.\ we assume that $\vek{x}(\bxi)$ and $\vek{f}(\bxi)$ are given by our stochastic
modelling.  As this are in total $M+N$ unknowns, we do not need more than
$L=M+N$ RVs $\bxi$.  The reason to not use $\vek{x}$ and $\vek{f}$ directly 
--- although that is certainly not excluded as e.g.\ for some $\ell$ the relation
$x_n=\xi_\ell$ is a possibility --- is that in the process of identification of
$\vek{x}$ or $\vek{f}$ they may turn out to be correlated, whereas the $\bxi$ can stay
independent as they are.  As now obviously
the state $\vek{u}(\bxi)$ is also a function of $\bxi$, the two \refeq{eq:model-N-p}
and \refeq{eq:resid-N-p} will then simply read
\begin{align} \label{eq:model-N} 
  \vek{A}(\vek{u}(\bxi),\vek{x}(\bxi)) &= \vek{f}(\bxi) ,\\
\label{eq:resid-N}   \vek{r}(\bxi) &= \vek{r}(\vek{u}(\bxi)) = 
             \vek{f}(\bxi) - \vek{A}(\vek{u}(\bxi),\vek{x}(\bxi)) = 0.
\end{align}
Computations such as e.g.\ evaluating the expected value of some function of the
response $\Psi(\vek{u},\vek{x})$ can then be transported to the variables
$\bxi=[\xi_1,\dots,\xi_\ell,\dots,\xi_L]^{\trpos}$
\begin{multline}  \label{eq:int-xi}    \EXP{\Psi(\vek{u},\vek{x})} = 
  \int_\Omega \Psi((\vek{u}(\omega),\vek{x}(\omega)))\, \D{P}(\di \omega) = 
  \int_\Omega \Psi((\vek{u}(\bxi(\omega)),\vek{x}(\bxi(\omega))))\,  \D{P}(\di \omega)\\ =
  \int_\Xi \Psi((\vek{u}(\bxi),\vek{x}(\bxi)))\, \Gamma(\di \bxi) =
  \int_\Xi \Psi(\bxi)\, \prod_{\ell=1}^L \Gamma_\ell(\di \xi_\ell)   \\=
  \idotsint  \Psi(\xi_1,\dots,\xi_\ell,\dots,\xi_L)\, 
      \Gamma_1(\di \xi_1)\dots\Gamma_\ell(\di \xi_\ell)\dots\Gamma_L(\di \xi_L),
\end{multline}
where the independence of the $\bxi$ allowed the use of Fubini's theorem to convert
the integral into a nested one-dimensional integration, and $\Gamma=\bxi_*\D{P}$ and the
$\Gamma_\ell = (\xi_\ell)_* \D{P}$ are the \emph{push-forward} or distribution
measures of the variables $\bxi$
and $\xi_\ell$, e.g.\ for normalised Gaussian variables $\Gamma_\ell(\di \xi_\ell)=
(2\uppi)^{-1/2}\exp(-\xi_\ell^2/2)\;\di \xi_\ell$ so that \refeq{eq:int-xi} can actually
be evaluated.

To actually describe the functions $\vek{u}(\bxi), \vek{x}(\bxi), \vek{f}(\bxi)$,
one further chooses a finite set of linearly independent functions 
$\{\psi_\alpha\}_{\alpha\in\C{J}_T}$ of the variables $\bxi(\omega)$, where
the index $\alpha=(\dots,\alpha_k,\dots)$ often is a \emph{multi-index},
and the set of multi-indices for approximation $\C{J}_T$ is a finite set with
cardinality (size) $T$.  Many different systems of functions can be used, classical
choices are \citep{wiener38, ghanemSpanos91, kallianpur, janson97, matthies6,
xiu2010numerical, knio2010spectral, Sullivan2015} multivariate polynomials ---
leading to the \emph{polynomial chaos expansion} (PCE) or generalised PCE (gPCE)
\citep{xiuKarniadakis02a}, as well as trigonometric functions \citep{janson97},
kernel functions as in kriging \citep{berlinet}, radial basis functions
\citep{Buhmann2000, Buhmann2003}, sigmoidal functions as in artificial neural
networks (ANNs) used in machine learning \citep{Murphy2012}, or functions derived
from fuzzy sets.  The particular choice is immaterial for the further development.
But to obtain results which match the above theory as regards $\E{L}$-invariant
subspaces, we shall assume that the set $\{\psi_\alpha\}_{\alpha\in\C{J}_T}$ includes
all the \emph{linear} functions of $\bxi$.  This is easy to achieve with
polynomials, and w.r.t\ kriging it corresponds to \emph{universal} kriging.
All other functions systems can also be augmented by a linear trend.

Thus a RV $\vek{u}(\bxi)$ would be replaced by a \emph{functional approximation}
--- this gives these methods its name, sometimes also termed \emph{spectral} approximation ---
\begin{equation}  \label{eq:FA}
   \vek{u}(\omega) = \sum_{\alpha\in\C{J}_T} \vek{u}_\alpha 
     \psi_\alpha(\bxi(\omega)) = \sum_{\alpha\in\C{J}_T} \vek{u}_\alpha 
     \psi_\alpha(\bxi) = \vek{u}(\bxi) .
\end{equation}

We describe the \emph{input} to the computational
model \refeq{eq:model-1}, namely $\vek{x}_f$, in a completely analogous way:
\begin{equation}  \label{eq:FA-x}
   \vek{x}_f(\omega) = \vek{x}_f(\bxi) = 
   \sum_{\alpha\in\C{J}_T} \vek{x}_\alpha \psi_\alpha(\bxi).
\end{equation}
Note that the parameters $\vek{x}$ are on purpose \emph{not} used as descriptive
variables, but rather $\bxi$, as later analysis may show that the
parameters $\vek{x}$ --- which are estimated --- are not independent.

The response of the system $\vek{y}(\bxi)=Y(\vek{x})=Y(\bxi)$ now has to be approximated by an
\emph{emulator}, \emph{proxy-} or \emph{meta-}model, or a \emph{response surface}.
Often this is achieved by approximating the whole state $\vek{u}(\bxi)$ by such
a proxy-model.  This is part of \emph{uncertainty quantification}
\citep{matthies6, xiu2010numerical, knio2010spectral, Sullivan2015}, and
is a computationally demanding task.  This produces
\begin{equation}  \label{eq:FA-y}
   Y(\vek{x}_f(\bxi)) = y(\bxi) = 
       \sum_{\beta\in\C{J}_K} \vek{y}_\beta \phi_\beta(\bxi),
\end{equation}
where of course the functions $\{\phi_\beta\}_{\beta\in\C{J}_K}$ could be
the same as $\{\psi_\alpha\}_{\alpha\in\C{J}_M}$, which we will assume from here on.
Similarly the error $\beps(\omega)$ has to be described, typically in RVs $\bbeta(\omega)=
[\eta_1(\omega),\dots,\eta_K(\omega)]$ independent of the RVs $\bxi(\omega)$,
\begin{equation}  \label{eq:FA-e}
   \beps(\bbeta) = \sum_{\gamma\in\C{J}_N} \beps_\gamma 
      \vphi_\gamma(\bbeta),
\end{equation}
where again the set of functions $\{\vphi_\gamma\}_{\gamma\in\C{J}_N}$ could be
the same as  $\{\phi_\beta\}_{\beta\in\C{J}_K}$ or $\{\psi_\alpha\}_{\alpha\in\C{J}_T}$.
In any case this gives
\begin{multline}  \label{eq:FA-z}
   \vek{z}(\bxi,\bbeta) = Y(\vek{x}(\bxi), \vek{u}(\bxi)) + \beps(\bbeta) =
   Y(\sum_{\alpha\in\C{J}_T} \vek{x}_\alpha \psi_\alpha(\bxi), 
   \sum_{\alpha\in\C{J}_T} \vek{u}_\alpha \psi_\alpha(\bxi)) + 
   \sum_{\gamma\in\C{J}_N} \beps_\gamma  \vphi_\gamma(\bbeta)\\
     = \vek{y}(\bxi)+\beps(\bbeta) = 
   \sum_{\beta\in\C{J}_K} \vek{y}_\beta \phi_\beta(\bxi) +
   \sum_{\gamma\in\C{J}_N} \beps_\gamma  \vphi_\gamma(\bbeta).
\end{multline}
As there is no loss in generality in assuming that all the functions are from the
same set, $\vphi_\alpha = \phi_\alpha = \psi_\alpha$, and a considerably simpler notation,
from now we shall do so.

In \refeq{eq:model-N} the space $\C{U}$ where the problem \refeq{eq:model-1}
for a fixed $\omega$ resp.\ $\bxi$ was formulated by an $N$-dimensional subspace
$\C{U}_N \subset \C{U}$.  Extending to the probabilistic description,
the solution $u(\omega)$ of \refeq{eq:model-1} lives in a tensor product space
$\C{U}\otimes\C{S}$, and the solution to \refeq{eq:model-N} hence lives in
the tensor product space $\C{U}_N\otimes\C{S}$.
By choosing a finite set $\{ \psi_\alpha \}_{\alpha\in\C{J}_T}$ of ansatz-functions
to represent all the RVs, one has defined an $T$-dimensional subspace
$\C{S}_T := \spn \{ \psi_\alpha \,:\, \alpha\in\C{J}_T\} \cong \D{R}^T$, $(T\in\D{N})$,
and the approximations mentioned above lie in the $(N\times T)$-dimensional subspace
\begin{equation}  \label{eq:fin-dim-space}
   \sum_{\alpha\in\C{J}_T} \vek{u}_\alpha \psi_\alpha(\bxi) \in \E{U}_{N,T} :=
   \C{U}_N\otimes\C{S}_T \subset \C{U}_N\otimes\C{S} \subset \C{U}\otimes\C{S} =: \E{U} .
\end{equation}

The RVs $\vek{x}(\bxi), \vek{f}(\bxi)$, and $\beps(\bbeta)$ are an input to the problem,
hence the coefficients in \refeq{eq:FA-x}, \refeq{eq:FA-e}, and similarly for the rhs
$\vek{f}(\bxi)$ can be considered given,  but the coefficients
$\vek{u}_\alpha$ in \refeq{eq:FA} or \refeq{eq:fin-dim-space} have to be computed.

\subsection{Intrusive or non-intrusive?}      \label{SS:intrusive}
Once one has decided what type of functions to use for approximation in the
proxy model, i.e.\ the subspace $\C{S}_T \subset \C{S}$ has been picked,
one has to decide how to determine the coefficients.  One of the
distinctions in the different methods is about what is to be evaluated.
One of the earliest and still
most frequent methods is to \emph{sample} the solution from \refeq{eq:model-N}
--- just like for the EnKF in \refeq{eq:EnKF} ---
at points $\bxi_s = \bxi(\omega_s) \in \Xi$ for
$\vek{u}(\bxi_s)$.  These are normal solves of \refeq{eq:model-N} for
certain realisations $\bxi_s \in \Xi$.
The points $\bxi_s$ may be chosen at random according to the measure
$\Gamma(\di \bxi)$ in \refeq{eq:int-xi} like in the
Monte Carlo method \citep{Hastings1970, schuellerSpanos01}, or according to
some deterministic quadrature rule like the quasi Monte Carlo method \citep{caflisch98}.

What is meant by the connotation in the title --- \emph{intrusive} or not? --- is that,
as is often the case, there
is software available to solve the deterministic problem in question, i.e.\ to compute
the solution $\vek{u}(\bxi_s)$ for a particular realisation $\bxi_s$.
Methods which only use this capability have then been termed ``non-intrusive'', as this
means that the underlying software does not have to be modified.
Without mentioning it there, it is obvious that the samplin methods described above
use only this capability, and are hence non-intrusive.  We shall see
that the methods to be described in \refSS{u-eval} also obviously fall into this
class.  It is also clear, as the computations for each realisation $\bxi_s$
can be performed independently, that the computation of all the values 
$\{ \vek{u}(\bxi_s) \}_{s=1}^S$ is ``embarassingly parallel''.
After this parallel phase, the results have to be summed or otherwise
post-processed in \refSS{u-eval}, and this phase can not be parallelised so easily.

Unfortunately, this denomination ``non-intrusive'' for the methods relying on the
evaluation of $\vek{u}(\bxi_s)$ through solution of \refeq{eq:model-N}
for realisations $\bxi_s$ could be understood to mean that other methods,
like the ones in \refSS{res-eval} which rely on evaluations of the
residuum $\vek{r}(\bxi_s)$ in \refeq{eq:resid-N}, are ``intrusive''
and actually \emph{do} require a modification of the underlying software.  This is
\emph{not} the case, as will be sketched later in \refSS{res-eval}.
But there is a distinction:  in the sampling methods above and
for the methods in \refSS{u-eval}
the evaluation of $\vek{u}(\bxi_s)$ for one particular realisation $\bxi_s$
is uncoupled from the evaluation for other realisations, hence the easy parallelism.
This is different in the methods in \refSS{res-eval}, here the evaluations
are coupled, one has to solve a \emph{coupled} system of equations.
Therefore the dichotomy should better be labelled ``uncoupled'' and ``coupled''.
In the literature for coupled systems (e.g.\ \cite{matthiesStf03}) the distinction is
the between ``monolithic'' methods, which indeed require typically modifications in
the software, and ``partitioned'' methods, which do not require modifications.
It is well-known that all coupled system may be also be solve in a partitioned way,
i.e.\ ``non-intrusively'' \cite{matthiesStf03}.

\subsection{Evaluating the solution for functional approximation}      \label{SS:u-eval}
Like always, there are several alternatives to determine the coefficients
$\vek{u}_\alpha$ in \refeq{eq:FA} for $u(\bxi)$,
in order to get a representation of the solution resp.\ state of the system.
Some of the possibilities when evaluating $\vek{u}(\bxi_s)$ are:
\begin{description}
\item[Stochastic collocation / interpolation:]
  This is one of the simplest ideas:  Compute $\vek{u}(\bxi_s)$ for particular
  points $\bxi_s \in \Xi$, and then interpolate with the functions $\psi_\alpha$.
  This is detailed in \refSSS{colloc}.
\item[Projection in function space:]
  This approach uses the idea to compute the coefficients $\vek{u}^{(\beta)}$
  through orthogonal projection in the Hilbert space $\C{S}$ onto the basis $\psi_\alpha$.
  This will be addressed in \refSSS{regress-project}.
\item[Discrete regression:]
  This is usually a combination of the interpolation and projection ideas.  Pure
  interpolation may suffer from so-called over-fitting, hence one uses more points
  $\btheta_s \in \Xi$ than there are functions $\psi_\alpha$, and then computes
  a least-squares approximation to the overdetermined system for the coefficients
  $\vek{u}^{(\beta)}$.  This is a discrete projection, and will be treated
  also in \refSSS{regress-project}.
\end{description}

\subsubsection{Stochastic collocation and interpolation}      \label{SSS:colloc}
As already stated, in this approach
\cite{babuska2007stochastic, xiu2005highorder, nobile2008sparse}
one computes the solution $\vek{u}(\bxi_s)$ in the interpolation point
$\{\bxi_s\}_{s=1}^S$.
These are deterministic solves at the sampling points $\bxi_s$.  Thus there is no
interaction of these solves, again one has only to solve many small systems
the size of the deterministic system.  This is very similar to the original
way of computing response surfaces.
The determining equations are
\begin{equation}  \forall s=1,\dots,S: \quad
\vek{u}(\bxi_s) = \sum_{\beta \in \C{J}_T} \vek{u}^{(\beta)} \psi_\beta(\bxi_s) . 
\label{eq:interp-cond}
\end{equation}
To write this in a concise form, the simplest is to consider this equation for
each component $u_n(\bxi)$ from $\vek{u}(\bxi) = [u_1(\bxi), \dots, u_n(\bxi),
\dots, u_N(\bxi)]^{\trpos}\in\D{R}^N$ and $u_n^{(\beta)}$ from the coefficient vector
$\vek{u}^{(\beta)} =[u_1^{(\beta)},\dots,u_n^{(\beta)},\dots,u_N^{(\beta)}]^{\trpos} 
\in \D{R}^N$; defining the $S\times T$ matrix 
\begin{equation}  \label{eq:interp-matx}
  \vek{\Psi}_{s,\beta} = \psi_\beta(\btheta_s),
\end{equation}
and the vectors $\vek{u}_n = [\dots,u_n^{(\beta)},\dots]_{\beta\in\C{J}_T}^{\trpos}
\in \D{R}^T$, and $\vek{y}_n = [\dots,u_n(\btheta_s),\dots]_{s=1,\dots,S}^{\trpos}\in\D{R}^S$,
the \refeq{eq:interp-cond} may be written as
\begin{equation}  \label{eq:interp-cond-m}  \forall n=1,\dots,N: \quad
   \vek{\Psi} \vek{u}_n = \vek{y}_n . 
\end{equation}
For the solution one requires that the matrix $\vek{\Psi}$ has to be non-singular.
The first condition for this is obviously that  $S=M$, and that
the functions $\{ \psi_\beta \}$ are linearly independent.
The second condition is that the points $\{ \btheta_s \}_{s=1}^S$
are \emph{uni-solvent} for the $\{\psi_\beta \}_{\beta\in\C{J}_T}$; this is equivalent with
the regularity of the matrix $\vek{\Psi}$.  If the system of points $\bxi_s$
and functions $\psi_\beta$ satisfies the so-called Kronecker-$\updelta$ condition,
\[ \vek{\Psi}_{s,\beta} = \psi_\beta(\btheta_s) = \updelta_{s,\beta}, \]
the system is particularly easy to solve, as $\vek{\Psi} = \vek{\Id}$, the identity
matrix.

One danger in interpolation is \emph{over-fitting}, where little errors in the evaluation
of  $\vek{u}(\btheta_s)$ are amplified.  Therefore one often uses more ``interpolation points''
than unknowns ($S>T$), leading to least-squares regression, see \refSSS{regress-project}.

\subsubsection{Spectral projection and regression}      \label{SSS:regress-project}
Another idea how to determine the coefficients $\vek{u}^{(\beta)}$ in
\refeq{eq:FA} is to project the solution $\vek{u}(\bxi)$ onto
the subspace $\spn\{ \psi_\beta \}_{\beta\in\C{J}_T} = \C{S}_T$, or rather
$\C{U}_N\otimes\C{S}_T$.  The simplest
way to achieve this to choose a set of linearly independent set of functions
$\{ \vphi_\alpha \}_{\alpha\in\C{J}_T}$, and to project along
$\spn\{ \vphi_\beta\, :\, \beta\in\C{J}_T \}$.  The orthogonality conditions are then
\begin{equation} \label{eq:orthog-proj-gen} \forall \alpha\in\C{J}_T: \quad
    \EXP{\vphi_\alpha(\cdot) \left( \vek{u}(\cdot) - 
    \sum_{\beta} \vek{u}^{(\beta)} \psi_\beta(\cdot) \right)} = 0
\end{equation}
Often one chooses an orthogonal projection by setting 
$\vphi_\alpha = \psi_\alpha$.  Then one has to solve the following system,
best again written for each component $u_n(\bxi)$ like in \refSSS{colloc}, using the
vector $\vek{u}_n$ from \refSSS{colloc}.  One additionally needs a new
rhs for each $n$, $\vek{v}_n = [\ldots, \EXP{\psi_\alpha(\cdot) u_n(\cdot)},
\ldots]_{\alpha\in\C{J}_T}^{\trpos}\in \D{R}^T$ and the $T\times T$ Gram matrix
$\vek{\Phi}_{\alpha,\beta} = \EXP{\psi_\alpha(\cdot)\psi_\beta(\cdot)}$:
\begin{equation} \label{eq:orthog-proj-cont} \forall n=1,\dots,N: \quad
    \vek{\Phi} \vek{u}_n = \vek{v}_n .
\end{equation}
This equation is equivalent to the minimising condition for the least-squares
solution of 
\[  \EXP{\nd{  \vek{u}(\cdot) - 
    \sum_{\beta} \vek{u}_n^{(\beta)} \psi_\beta(\cdot)}^2} \to \min ,
\]
defining an orthogonal projection in $\C{S}$ onto $\C{S}_T$.
If one then uses a numerical quadrature rule to compute the Gram matrix $\vek{\Phi}$,
\begin{equation} \label{eq:num-Gram} 
    \vek{\Phi}_{\alpha,\beta} = \EXP{\psi_\alpha \psi_\beta} = 
    \int_{\Xi} \psi_\alpha(\bxi) \psi_\beta(\bxi)\; \Gamma(\di \bxi)
    \approx \sum_{s=1}^S w_s^2 \psi_\alpha(\bxi_s) \psi_\beta(\bxi_s)
\end{equation}
with sampling points $\bxi_s\in\Xi$ and positive weights $w_s^2$, the numerical quadrature
equivalent of \refeq{eq:orthog-proj-cont} is
\begin{equation} \label{eq:orthog-proj-discr} \forall n=1,\dots,N: \quad
    \vek{\Psi}^{\trpos} \vek{W}^2 \vek{\Psi} \vek{u}_n = \vek{\Psi}^{\trpos} \vek{W}^2 \vek{y}_n ,
\end{equation}
where again $(\vek{\Psi}_{s,\alpha}) = (\psi_\alpha(\bxi_s))^{\trpos} \in \D{R}^{S\times T}$,
$\vek{W}=\diag(w_s)\in\D{R}^{S\times S}$, and
\[ \vek{\Phi} \approx \vek{\Psi}^{\trpos} \vek{W}^2 \vek{\Psi}, \text{ and  } 
 \vek{y}_n = [\ldots,u_n(\bxi_s),\ldots]^{\trpos}_{s=1,\dots,S} . \]
The rhs $\vek{\Psi}^{\trpos} \vek{W}^2 \vek{y}_n$ comes from
\[ \forall \alpha\in\C{J}_T: \quad
  (\vek{v}_n)_\alpha = \EXP{\psi_\alpha(\cdot) u_n(\cdot)} \approx 
     \sum_{s=1}^S w_s^2 \psi_\alpha(\bxi_s) u_n(\bxi_s) = 
     (\vek{\Psi}^{\trpos} \vek{W}^2 \vek{y}_n)_\alpha .
\]

A similar idea, a projection, is typically used in the case that one uses more
points $\{ \bxi_s \}_{s=1}^S$ than functions $\{\psi_\beta\}_{\beta\in\C{J}_T}$
in \refSSS{colloc}.  A least-squares approach to \refeq{eq:interp-cond-m} then yields
\begin{equation} \label{eq:regr-lsqr}  \forall n=1,\dots,N: \quad
  \nd{\vek{\Psi} \vek{u}_n - \vek{y}_n}_{\ell_2^S(\vek{W})^2}^2 \to \min.
\end{equation}
Here the discrete norm $\nd{\cdot}_{\ell_2^S(\vek{W}^2)}$ is with weights 
$\vek{W}^2=\diag(w_s^2)$, so that $\nd{\vek{v}_n}_{\ell_2^S(\vek{W}^2)}^2 =
\sum_s w_s^2 (v^{(s)}_n)^2$.  The Galerkin condition of \refeq{eq:regr-lsqr} is then
\begin{equation} \label{eq:orthog-proj-lsqr} \forall n: \quad
    \vek{\Psi}^{\trpos} \vek{W}^2 \vek{\Psi} \vek{u}_n = \vek{\Psi}^{\trpos} \vek{W}^2 \vek{y}_n ,
\end{equation}
One may observe the direct similarity of \refeq{eq:orthog-proj-lsqr} with
\refeq{eq:orthog-proj-discr}.  Hence the discrete least-squares regression
\refeq{eq:regr-lsqr} with discrete weights $w_s^2$
may be interpreted, in case the interpolation points $\bxi_s$ are
quadrature points of a numerical integration rule with weights $w_s^2$,
as a quadrature approximation of the continuous case \refeq{eq:orthog-proj-cont}.

Numerically it may well happen that the least-squares matrix
$\vek{\Psi}^{\trpos} \vek{W}^2 \vek{\Psi}$
in either \refeq{eq:orthog-proj-discr} or \refeq{eq:orthog-proj-lsqr} is ill-conditioned.
Then it may be more advisable \cite{Golub1996}, instead of solving the systems
by e.g.\ Choleski-factorisation of $\vek{\Psi}^{\trpos} \vek{W}^2 \vek{\Psi}$, to compute
the least-squares solution of the the ``square-root'' systems
\begin{equation} \label{eq:pro-reg-sqr-rt} \forall n: \quad
    \vek{W} \vek{\Psi} \vek{u}_n =  \vek{W} \vek{y}_n 
\end{equation}
through e.g.\ QR-decomposition of the matrix $\vek{W} \vek{\Psi}$.  The condition number
of $\vek{W} \vek{\Psi}$ is only the square root of 
the condition number of $\vek{\Psi}^{\trpos} \vek{W}^2 \vek{\Psi}$.

Hence both in this section as well as in \refSSS{colloc}, the coefficients 
$\vek{u}^{(\beta)}$ in \refeq{eq:FA} are essentially computed by evaluating
the solution $\vek{u}(\bxi_s)$ at discrete points $\bxi_s$.

\subsection{Galerkin methods for functional approximation}        \label{SS:res-eval}
This is the the approach to use $\vek{r}(\bxi)$ from
\refeq{eq:resid-N} to evaluate the coefficients.  Actually, that equation resulted
for a fixed $\bxi$ resp. $\omega$ from  projecting \refeq{eq:model-1} in $\C{U}$
onto the $N$-dimensional subspace $\C{U}_N$.  Here we project in some ways
additionally onto $\C{S}_T$, hence in total onto $\E{U}_{N,T} = \C{U}_N\otimes\C{S}_T$.
\begin{description}
\item[Least-Squares:]
  Here the \refeq{eq:resid-N} is considered as an element $\vek{r}(\vek{u}(\bxi))$
  in the Hilbert space $\C{U}_N^*\otimes\C{S}_T$, which should vanish at the solution
  $\vek{u}(\bxi)=\sum_\beta \vek{u}^{(\beta)}\psi_\beta(\bxi)\in \E{U}_{N,T}$.
  So one can compute the norm squared of the residuum, and then choose
  the coefficients $\vek{u}^{(\beta)}$ in such a way as to minimise this.
\item[Galerkin-Methods:]
  The Galerkin idea is to project the  residuum $\vek{r}(\bxi)$ onto the
  subspace $\C{U}_N\otimes \C{S}_T$.  This gives the condition
  --- Galerkin orthogonality --- to determine the coefficients $\vek{u}^{(\beta)}$.
  Most often one uses an orthogonal projection.
  This is considered further here.
\end{description}

These groups of methods, as well least-squares methods, start from
\refeq{eq:resid-N}.  The least squares solution  will not be considered
further, but we concentrate on Galerkin methods.  Similarly to
\refSSS{regress-project}, a projection is computed.  Here it is not the
solution which is projected, but the residuum.  As in \refSSS{regress-project},
one chooses a set of linearly independent set of functions
$\{ \vphi_\alpha \}_{\alpha\in\C{J}_T}$ to project along their span.

The Galerkin condition of \refeq{eq:resid-N} is then
\begin{equation}  \label{eq:srhs-res-G} \forall \alpha\in\C{J}_T: \quad
 \EXP{\vphi_\alpha(\cdot)\vek{r}(\cdot)} = \EXP{\vphi_\alpha(\cdot) \left(
 \vek{f}(\cdot) - \vek{A}(\cdot,\sum_{\beta\in\C{J}_M}\vek{u}^{(\beta)}\psi_\beta(\cdot))\right)}
  = 0 . 
\end{equation}
Often one wants an orthogonal projection by setting $\vphi_\alpha = \psi_\alpha$.
In any case, the \refeq{eq:srhs-res-G} defines a system of $T$ equations 
of dimension $N$ to determine the $T$ coefficients $\vek{u}^{(\beta)}\in\C{U}_N$.

The \refeq{eq:resid-N} results in the space $\E{U}_{N,T}$ in
\begin{equation}  \label{eq:res_nonlin_p}
  \tnb{r}(\tnb{u}) = [\ldots, \EXP{\psi_\alpha(\cdot)
  \vek{r}(\cdot)}, \ldots]_{\alpha\in\C{J}_T} = 0,
\end{equation}
where the same block vectors --- $\tnb{u} = [\ldots, \vek{u}^{(\alpha)},
\ldots]_{\alpha\in\C{J}_T}$ --- as before are used.
Now \refeq{eq:res_nonlin_p} is a huge non-linear system of dimension $N\times T$,
and one way to approach it is through the use of Newton's method, which involves
linearisation and subsequent solution of the linearised system.  Differentiating
the residuum in \refeq{eq:res_nonlin_p}, one obtains for the
$(\alpha, \beta)$ block-element of the derivative
\begin{equation}  \label{eq:res_nonlin_deriv}
  (\Di \tnb{r}(\tnb{u}))_{\alpha \beta} = \EXP{\psi_\alpha \, [\Di
  \vek{r}(\vek{u}(\bxi))]\, \psi_\beta} = - \EXP{\psi_\alpha \, [\Di_u
  \vek{A}(\bxi,\vek{u}(\bxi))]\, \psi_\beta}.
\end{equation}
Denote the matrix with the entries in \refeq{eq:res_nonlin_deriv} by
$-\tnb{K}_T(\tnb{u})$.  It is a tangential stiffness matrix.  If we are to use
Newton's method to solve the nonlinear system \refeq{eq:res_nonlin_p}, at
iteration $k$ it would look like
\begin{equation}  \label{eq:NR-it}
  \tnb{K}_T(\tnb{u}_k) (\tnb{u}_{k+1} - \tnb{u}_k) = \tnb{r}(\tnb{u}_k).
\end{equation}
One may now use all the techniques developed for linear problems so far to solve this
equation, and this then really is the workhorse for the non-linear equation.

Another possibility, avoiding the costly linearisation and solution of a new
linear system at each iteration, is the use of limited memory quasi-Newton methods
\cite{matthiesStrang79}.  This was done in \cite{keeseMatthies03siam}, and the
quasi-Newton method used---as we have a symmetric positive definite or
potential minimisation problem this was the {\em BFGS}-update---performed very
well.  The quasi-Newton methods produce updates to the inverse of a matrix,
and these low-rank changes are also best kept in tensor product form
\cite{matthiesStrang79}; so that we have tensor products here on two levels,
which makes for a very economical representation.

But in any case, in each iteration the residual \refeq{eq:res_nonlin_p} has
to be evaluated at least once, which means that for all $\alpha\in{\C J}_{T}$
the integral
\begin{equation}  \label{eq:resid-nonl-ex}
\EXP{\psi_\alpha(\cdot) \vek{r}(\cdot)} = \int_{\Xi}
\psi_\alpha(\bxi) \vek{r}(\bxi)\, \Gamma(\di \bxi)
\end{equation}
has to be computed.  In general this can not be done analytically,
and one has to resort to numerical quadrature rules:
\begin{equation}  \label{eq:resid-nonl-apprx}
\int_{\Xi} \psi_\alpha(\bxi) \vek{r}(\bxi)\,
\Gamma (\di \bxi) \approx \sum_{s=1}^S w_s \psi_{\alpha}(\bxi_s) \vek{r}(\bxi_s) .
\end{equation}
What this means is that for each evaluation of the residual
\refeq{eq:res_nonlin_p} the spatial residuum \refeq{eq:resid-N} has to be
evaluated $S$ times --- once for each $\bxi_s$ where one has to compute
$\vek{r}(\bxi_s)$.  Certainly this can be done
independently and in parallel without any communication.  We additionally
would like to point out that instead of solving the system every time for each
$\bxi_s$ as in the methods in \refSS{u-eval}, here one
only has to compute the residuum --- in fact typically a preconditioned 
residuum by performing one iteration --- at $\bxi_s$, which is typically
much cheaper.  This formulation with numerical integration makes the
Galerkin method completely ``non-intrusive'' \cite{LgAlDlHgmAn13-p}.
In fact this is a partitioned solution of the coupled set of equations
\refeq{eq:res_nonlin_p}.  This can be further extended also to compute
low-rank approximations to the solution directly in a non-intrusive way;
in \cite{LgDlHgmAn14-p} this is shown for the proper generalised decomposition
(PGD) in conjunction with BFGS iterations.

\subsection{The functional or spectral Kalman filter}        \label{SS:SPKF}
We come back to the task of computing the filter \refeq{eq:GMKF-RV}, where
we inject the terms from \refeq{eq:FA-z}, so that with the now hopefully
determined coefficients $\vek{u}_\alpha$ in \refeq{eq:FA} of the solution and
\[ 
  z(\bxi,\bbeta) = Y(\sum_{\alpha\in\C{J}_T} \vek{x}_{f,\alpha} \psi_\alpha(\bxi),
    \sum_{\alpha\in\C{J}_T} \vek{u}_\alpha \psi_\alpha(\bxi)) +
  \sum_{\gamma\in\C{J}_T} \beps_\gamma  \vphi_\gamma(\bbeta)  = 
  \sum_{\alpha\in\C{J}_T} \vek{y}_\alpha \psi_\alpha(\bxi) +
  \sum_{\gamma\in\C{J}_T} \beps_\gamma  \vphi_\gamma(\bbeta)
\]
it reads 
\begin{multline} \label{eq:PCEKF}
  \vek{x}_a(\bxi,\bbeta) =  \vek{x}_f(\bxi) + 
      \vek{C}_{x_f z} \vek{C}^{-1}_{z}(\vcek{y}-z(\bxi,\bbeta))
      =  \vek{x}_f(\bxi) + \vek{K}(\vcek{y}-z(\bxi,\bbeta)) \\
      = \sum_{\alpha\in\C{J}_T} \vek{x}_{f,\alpha} \psi_\alpha(\bxi)
      + \vek{K}\left(\vcek{y}- \left(\sum_{\alpha\in\C{J}_T} \vek{y}_\alpha \psi_\beta(\bxi)
      + \sum_{\gamma\in\C{J}_T} \beps_\gamma  \vphi_\gamma(\bbeta) \right)\right).
\end{multline}
This has been termed --- especially if the approximating functions are polynomials ---
as a \emph{polynomial chaos expansion Kalman filter}; a better name is the
\emph{spectral Kalman filter} (SPKF).  This is
an explicit and easy to evaluate expression for the assimilated or
\emph{updated} variable in terms of the input and the state $\vek{u}(\bxi)$.

It remains to show how to approximate the Kalman gain operator $\vek{K}=
\vek{C}_{x_f z} \vek{C}^{-1}_{z}$.  This actually
is fairly straightforward with functional approximations.  One has
\begin{equation}  \label{eq:Cxz-FA}
   \vek{C}_{x_f z} = \EXP{\vtil{x}_f(\cdot) \otimes \vtil{z}(\cdot,\cdot)}
   = \EXP{\vtil{x}_f(\cdot)\otimes(\vtil{y}(\cdot) + \beps(\cdot))}=
   \vek{C}_{x_f y} +    \vek{C}_{x_f \vepsilon},
\end{equation}
where the last term will vanish if $\bxi$ and $\bbeta$ are independent.
Further one has
\begin{equation}  \label{eq:Czz-FA}
   \vek{C}_{z} = \EXP{\vtil{z}(\cdot,\cdot) \otimes \vtil{z}(\cdot,\cdot)}
   = \EXP{(\vtil{y}(\cdot) + \beps(\cdot)) \otimes (\vtil{y}(\cdot) + \beps(\cdot))} =
     \vek{C}_{y} + \vek{C}_{\vepsilon y} + \vek{C}_{\vepsilon y}^T + \vek{C}_{\vepsilon} ,
\end{equation}
where again the two middle terms will vanish if $\bxi$ and $\bbeta$ are independent.
Assuming this, one has
\begin{align}  \label{eq:Cxx-FA}
  \vek{C}_{x_f y}&=\sum_{\alpha,\beta \in \C{J}_T} \EXP{\psi_\alpha(\cdot)\psi_\beta(\cdot)}
         \vek{x}_{f \alpha} \otimes \vek{y}_\beta  - \vbar{x} \otimes \vbar{y} ,\\
  \vek{C}_{\vepsilon}&=\sum_{\alpha,\beta \in \C{J}_T} 
     \EXP{\vphi_\alpha(\cdot)\vphi_\beta(\cdot)} 
     \beps_\alpha \otimes \beps_\beta  - \vbar{\vepsilon} \otimes \vbar{\vepsilon}, \\
  \vek{C}_{y}&=\sum_{\alpha,\beta \in \C{J}_T} \EXP{\psi_\alpha(\cdot)\psi_\beta(\cdot)}
      \vek{y}_\alpha \otimes \vek{y}_\beta  - \vbar{y} \otimes \vbar{y} .
\end{align}
Now all ingredients for the SPKF in \refeq{eq:PCEKF} are given explicitly in terms
of known coefficients and known RVs, and hence may be directly computed, and one has
an explicit expression for $\vek{x}_a(\bxi,\bbeta)$.

\subsection{Examples with the linear spectral filter}        \label{SS:ex-SPKF}
This is to show some examples computed with \refeq{eq:PCEKF}.
\begin{figure}[!ht]
\centering
 \includegraphics[width=0.8\textwidth,height=0.35\textheight]{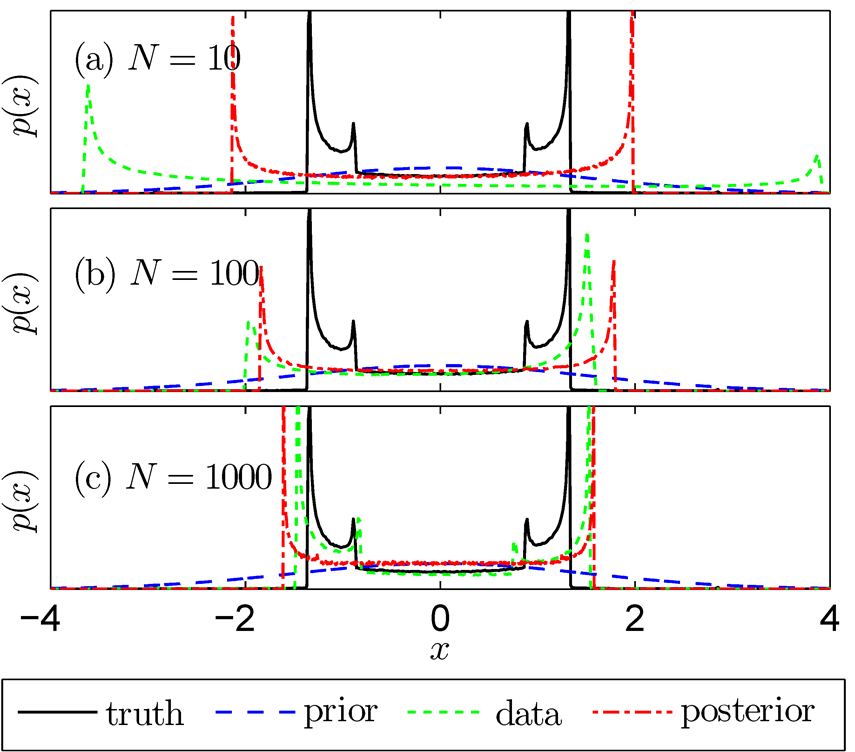}
 \caption{pdfs for linear Bayesian update, from \citep{opBvrAlHgm12}}
\label{F:exp-B-1}
\end{figure}
As the traditional Kalman filter is highly geared towards Gaussian
distributions \citep{Kalman}, and also its Monte Carlo variant EnKF
which was mentioned previously at the beginning of this section tilts towards Gaussianity,
we start with a case --- already described in \citep{opBvrAlHgm12} --- where
the the quantity to be identified has a strongly
non-Gaussian distribution, shown in black --- the `truth' --- in \refig{exp-B-1}.
The operator describing the system is the identity --- we compute the quantity
directly, but there is a Gaussian measurement error.  The `truth' was
represented as a $12^{\text{th}}$ degree PCE.

We use the methods as described in \refSS{SPKF}, and here in particular
the \refeq{eq:PCEKF}, the SPKF.
The update is repeated several times (here ten times) with new
measure\-ments---see \refig{exp-B-1}.  The task is here to identify the
distribution labelled as `truth' with ten updates of $N$ samples
(where $N=10, 100, 1000$ was used), and we start with a very broad
Gaussian prior (in blue).  Here we see the ability of the polynomial
based LBU, the PCEKF, to identify highly non-Gaussian distributions,
the posterior is shown in red and the pdf estimated from the samples
in green; for further details see \citep{opBvrAlHgm12}.

\begin{figure}[!ht]
\centering
 \includegraphics[width=0.9\textwidth,height=0.35\textheight]{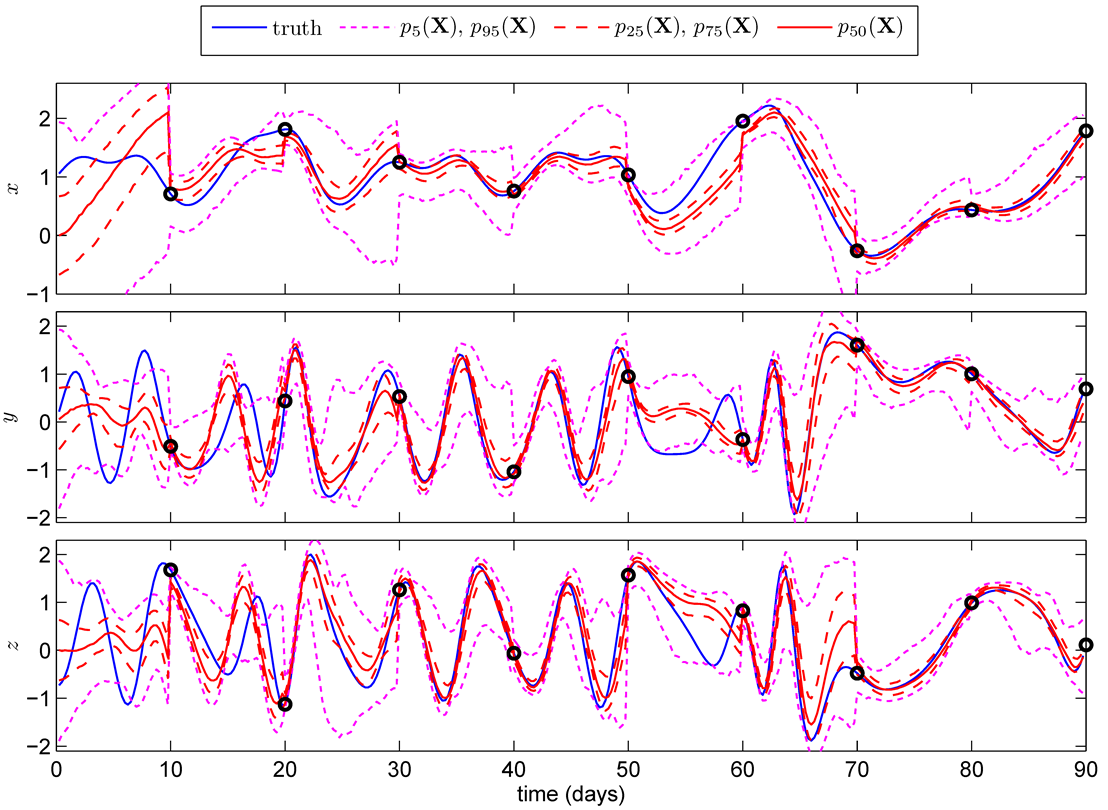}
 \caption{Time evolution of Lorenz-84 state and uncertainty with the LBU, from \citep{opBvrAlHgm12}}
\label{F:exp-B-2}
\end{figure}
The next example is also from \citep{opBvrAlHgm12}, where the system
is the well-known Lorenz-84 chaotic model, a system of three nonlinear
ordinary differential equations operating in the chaotic regime. This
is truly an example.  Remember that this was
originally a model to describe the evolution of some amplitudes of a
spherical harmonic expansion of variables describing world climate.
As the original scaling of the variables has been kept, the time axis
in \refig{exp-B-2} is in \emph{days}.  Every ten days a noisy
measurement is performed and the state description is updated.  In
between the state description evolves according to the chaotic dynamic
of the system.  One may observe from \refig{exp-B-2} how the
uncertainty --- the width of the distribution as given by the quantile
lines---shrinks every time a measurement is performed, and then
increases again due to the chaotic and hence noisy dynamics.  Of
course, we did not really measure world climate, but rather simulated
the `truth' as well, i.e.\ a \emph{virtual} experiment, like the
others to follow.  More details may be found in \citep{opBvrAlHgm12}
and the references therein.

\begin{figure}[!ht]
\centering
\begin{minipage}{.4\textwidth}
  \centering
  \includegraphics[width=.99\linewidth]{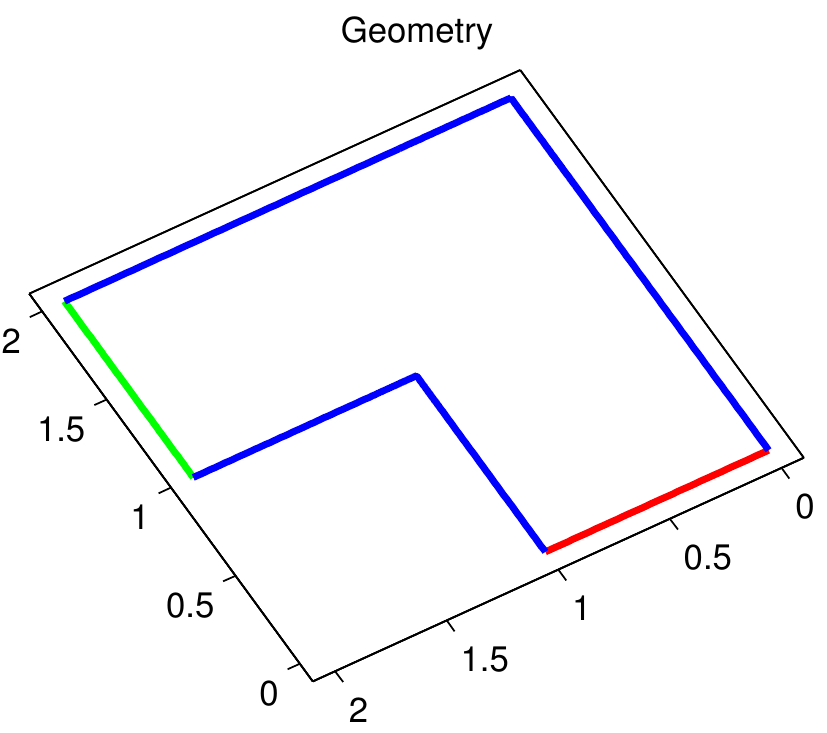}
  \caption{Diffusion domain, from \citep{bvrAlOpHgm12-a}}
  \label{F:exp-B-3}
\end{minipage}%
\begin{minipage}{.6\textwidth}
  \centering
  \includegraphics[width=.99\linewidth]{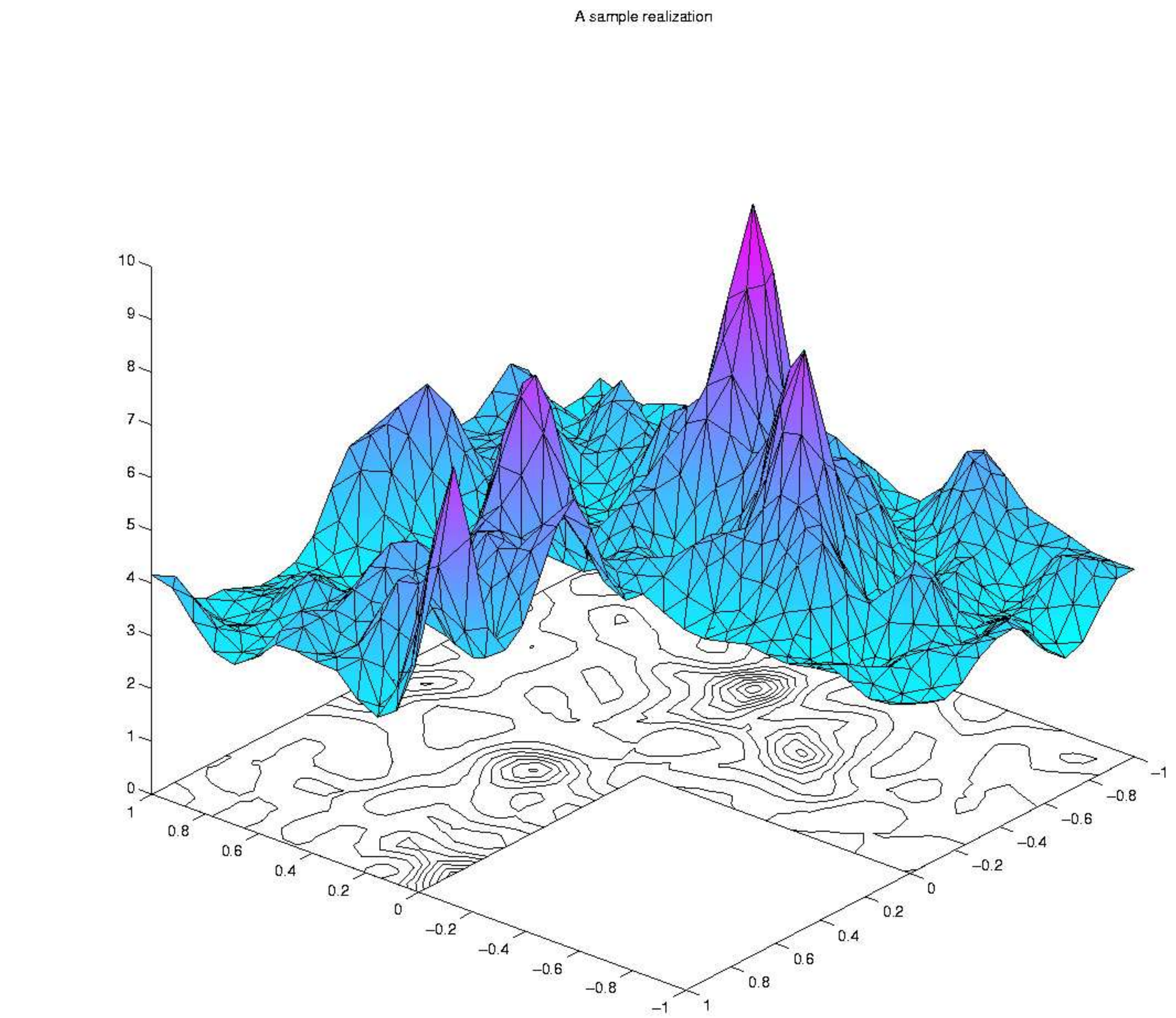}
  \caption{Conductivity field, from \citep{bvrAlOpHgm12-a}}  
  \label{F:exp-B-4}
\end{minipage}
\end{figure}
From \citep{bvrAlOpHgm12-a} we take the example shown in
\refig{exp-B-3}, a linear stationary diffusion equation on an L-shaped
plane domain.  The diffusion coefficient
$\kappa$  is to be identified.  As argued in \citep{BvrAkJsOpHgm11},
it is better to work with $q = \log \kappa$ as the diffusion coefficient has
to be positive, but the results are shown in terms of $\kappa$.

One possible realisation of the diffusion coefficient
is shown in \refig{exp-B-4}.  More realistically, one should assume that
$\kappa$ is a symmetric positive definite tensor field, unless one knows that
the diffusion is \emph{isotropic}.  Also in this case one should do the updating
on the logarithm.  For the sake of simplicity we stay with the scalar case,
as there is no principal novelty in the non-isotropic case.
The virtual experiments use different right-hand-sides $f$ in \refeq{eq:model-1},
and the measurement is the observation of the solution $u$ averaged over little patches.

\begin{figure}[!ht]
\centering
\begin{minipage}{0.48\textwidth}
  \centering
  \includegraphics[width=0.99\linewidth]{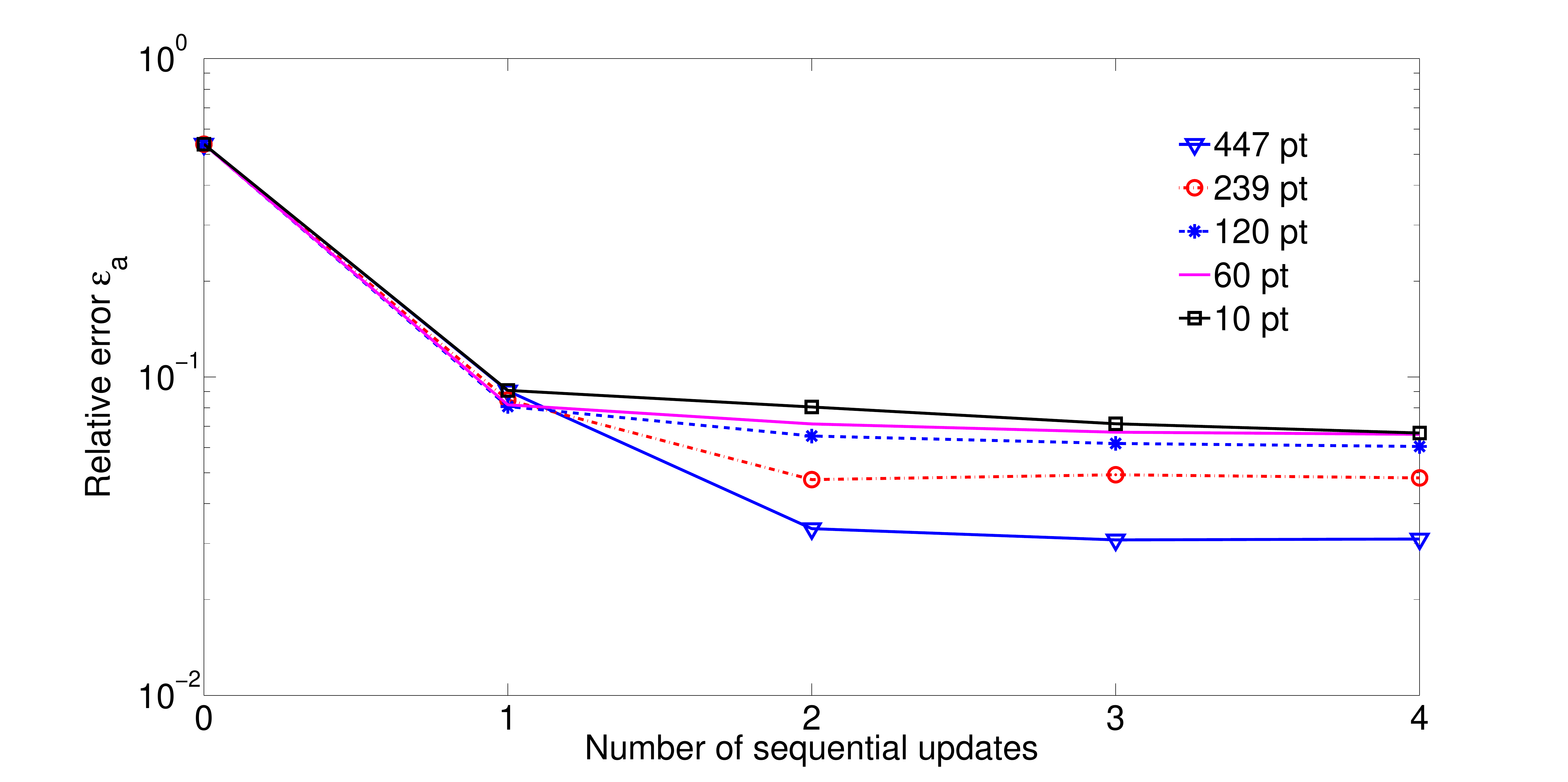}
  \caption{Convergence of identification, from \citep{bvrAlOpHgm12-a}}
  \label{F:exp-B-7}
\end{minipage}%
\hfill
\begin{minipage}{0.48\textwidth}
  \centering
  \includegraphics[width=0.99\linewidth,height=0.17\textheight]{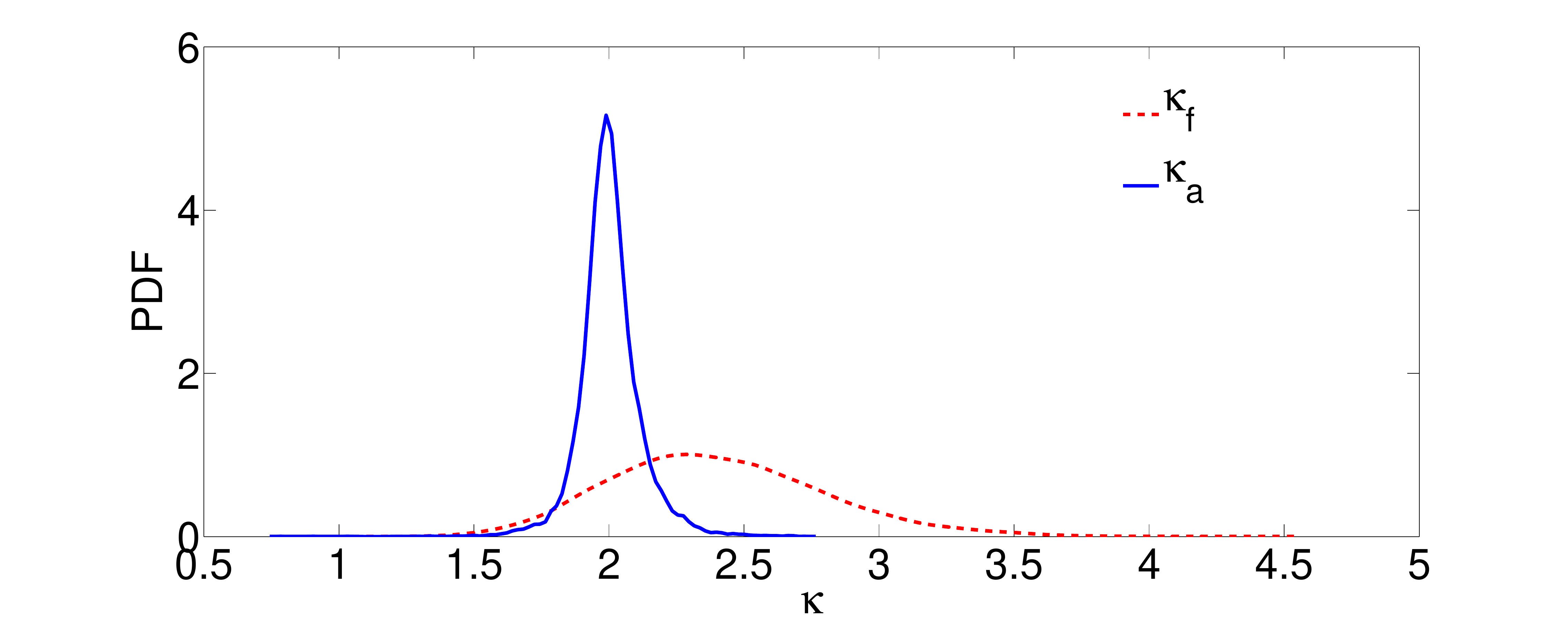}
  \caption{Prior and posterior, from \citep{bvrAlOpHgm12-a}}  
  \label{F:exp-B-8} 
\end{minipage}
\end{figure}
In \refig{exp-B-7} one may observe the decrease of the error with successive
updates, but due to measurement error and insufficient information from just
a few patches, the curves level off, leaving some residual uncertainty.
The pdfs of the diffusion coefficient at some point in the domain before
and after the updating is shown in \refig{exp-B-8}, the `true' value at
that point was $\kappa=2$.
Further details can be found in \citep{bvrAlOpHgm12-a}.

\section{More accurate filters}  \label{S:more-acc}
The  filter in \refeq{eq:GMKF} resp.\ \refeq{eq:PCEKF}
is in some ways the simplest possible one.
More accurate filters than the GMKF in \refeq{eq:GMKF} resp.\ \refeq{eq:GMKF-RV}
can be achieved in two ways:  for one the subspace $\E{X}_1\subset \E{X}_{\sigma(z)}$
in \refeq{eq:X-1} of affine maps of $\vek{z}$ may be replaced by increasingly
larger subspaces $\E{X}_\ell$, and hence more accurate optimal maps
$\vpi_{\Psi,\ell}(\vek{z})$ approximating $\vpi_\Psi(\vek{z})$.  This makes for
a better approximation of the conditional mean 
$\vbar{x}^{\vcek{y}}_{a,\ell} := \vpi_{x,\ell}(\vcek{y})$ to $\vbar{x}^{\vcek{y}}_a$.
On the other hand, if one wants $\tnb{x}_a$ to give a better approximation to
the Bayesian posterior, the zero-mean part $\vtil{x}^{\vcek{y}}_{a,\ell}$
has to be possibly transformed.

\subsection{Better approximation to the conditional expectation}  \label{SS:better-CE}
Let us start by choosing approximating subspaces $\E{X}_\ell \subset \E{X}_{\sigma(z)}$ with
\[  \E{X}_1\subset \E{X}_\ell \subset \E{X}_{\sigma(z)} \subset \E{X} . \]
For a RV $\Psi(\vek(x))$, this should give a better approximation $\vpi_{\Psi,\ell}(\vek{z})$
to $\vpi_\Psi(\vek{z})$ in \refeq{eq:opt-map} than the linear map in \refeq{eq:orth-CE-lin}.
Assuming that the subspaces $\E{X}_\ell$ are chosen such that their union is dense
in $\E{X}_{\sigma(z)}$,
\begin{equation} \label{eq:X-dense}
 \cl\left(\bigcup_{\ell=1}^\infty \E{X}_\ell\right) =  \E{X}_{\sigma(z)}, 
\end{equation}
one may approximate with $\vpi_{\Psi,\ell}$ the optimal map $\vpi_\Psi$ to
any desired accuracy by taking $\ell$ large enough.  This is shown
in  \citep{HgmEzBvrAlOp15, HgmEzBvrAlOp15-p} in general, and in particular
for the case when $\E{X}_\ell$ is the subspace given by polynomials up to degree $\ell$
in $\vek{z}$.

Using this in the case $\Psi(\vek{x}) = \vek{x}$,
the linear filter \refeq{eq:GMKF-RV} would then be replaced by
\begin{equation} \label{eq:CE-l-filter}
  \vek{x}_{a,\ell}(\omega) =  \vek{x}_f(\omega) + 
      \vpi_{x,\ell}(\vcek{y}) - \vpi_{x,\ell}(\vek{z}(\omega)) =
      \vpi_{x,\ell}(\vcek{y}) +(\vek{x}_f(\omega) - \vpi_{x,\ell}(\vek{z}(\omega))) ,
\end{equation}
which is a non-linear filter
as an approximation to the CEM-filter \refeq{eq:orth-CE-map}.  Observe that in
general this will only result in the RV $\vek{x}_{a,\ell}$ having a posterior
mean $\vpi_{x,\ell}(\vcek{y})=\vbar{x}^{\vcek{y}}_{a,\ell}$ closer to the posterior
Bayesian mean $\vbar{x}^{\vcek{y}}_a$.  In case the density condition \refeq{eq:X-dense}
is satisfied, one obtains convergence $\vbar{x}^{\vcek{y}}_{a,\ell}\to\vbar{x}^{\vcek{y}}_a$
as $\ell\to\infty$.

To introduce the nonlinear filter as just sketched, one may
look shortly at a very simplified example, identifying a value, where
only the third power of the value plus a Gaussian error RV is observed.
All updates follow \refeq{eq:CE-l-filter}, but the update map is computed
with different accuracy.
\begin{figure}[!ht]
\centering
 \includegraphics[width=0.9\textwidth,height=0.35\textheight]{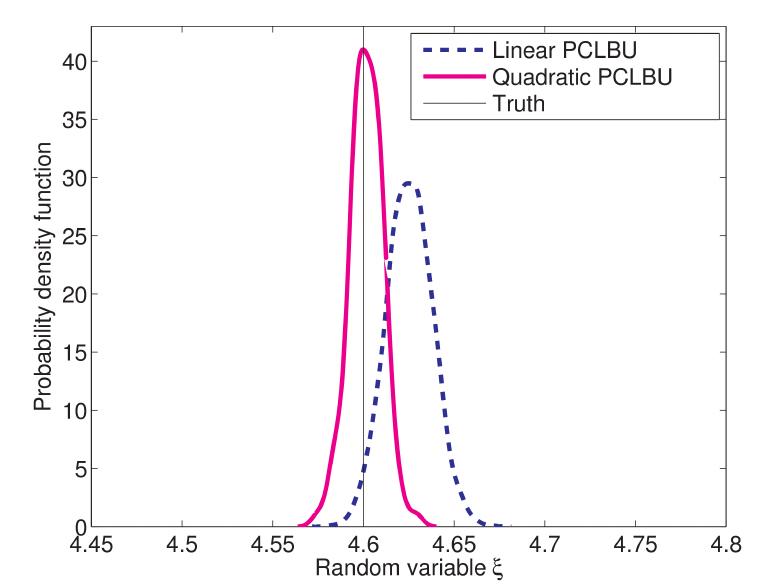}
 \caption[Updates for cubic observations]{Perturbed observations of the cube of a RV, different
 updates --- linear and quadratic update}
\label{F:exp-B-2}
\end{figure}
Shown are the pdfs produced by the linear filter according to \refeq{eq:PCEKF}
--- Linear polynomial chaos Bayesian update (Linear PCBU) ---
a special form of \refeq{eq:CE-l-filter}, 
and using polynomials up to order two, the quadratic polynomial chaos
Bayesian update (QPCBU).  One may observe that due to the nonlinear observation,
the differences between the linear filters and the quadratic one are already significant,
the QPCBU yielding a better update.

\subsection{Transformation of the zero-mean part}  \label{SS:trfm-zm}
Hence, if on the other hand one wants to construct a RV which matches the full posterior
Bayesian distribution, one has to look at the zero mean part from \refeq{eq:orth-CE-map}
\begin{equation} \label{eq:CE-zm}  
   \vtil{x}^{\vcek{y}}_{a}(\bxi,\bbeta) =  \vek{x}_f(\bxi) - \vpi_{x}(\vek{z}(\bxi,\bbeta)) ,
\end{equation}
which is essentially what is left from the orthogonal projection.
For an actual computation, one would choose a finite $\ell$ and use \refeq{eq:CE-l-filter}.
This RV $\vtil{x}^{\vcek{y}}_{a}$ will have to be transformed further.

From \refeq{eq:GMKF-var} one has the covariance $\vek{C}_{x_a}$ of
$\vtil{x}^{\vcek{y}}_{a,1}$.  A similar computation can be performed,
at least numerically, for $\vtil{x}^{\vcek{y}}_{a}$ in \refeq{eq:CE-zm},
giving $\vhat{C}_{x_a} =  \EXP{\vtil{x}^{\vcek{y}}_{a} \otimes \vtil{x}^{\vcek{y}}_{a}}$.

On the other hand, one may compute the correct value of the posterior covariance
to any degree of accuracy with an optimal map.  While in \refSSS{CE-filt} we
have computed the optimal map $\vpi_\Psi$ for $\Psi(\vek{x}) = \vek{x}$,
now we may take $\Psi(x)=\vek{x}\otimes\vek{x}$ to obtain an optimal map
$\vpi_{x\otimes x}$.  This then gives
\begin{equation} \label{eq:CE-covar-a}  
   \vek{C}_{x_a} = \EXP{\vek{x}_f \otimes \vek{x}_f | \vcek{y}} - 
   \vbar{x}^{\vcek{y}}_a \otimes \vbar{x}^{\vcek{y}}_a =
   \vpi_{x\otimes x}(\vcek{y}) - \vbar{x}^{\vcek{y}}_a \otimes \vbar{x}^{\vcek{y}}_a .
\end{equation}
For a numerical approximation over $\E{X}_\ell$ this would similarly result in
an approximation $\vek{C}_{x_a,\ell}$.
Both the former $\vhat{C}_{x_a}$ as well as $\vek{C}_{x_a}$ in \refeq{eq:CE-covar}
are spd matrices, hence have a square root.  Therefore, when considering the following RV
\begin{equation} \label{eq:CE-l-m-cov}
   \vek{x}_{a,\text{cov}}(\bxi,\bbeta) = \vpi_{x}(\vcek{y}) + 
   \vek{C}_{x_a}^{1/2} \,\vhat{C}_{x_a}^{-1/2}(\vtil{x}^{\vcek{y}}_{a}(\bxi,\bbeta)) ,
\end{equation}
with $\vtil{x}^{\vcek{y}}_{a}$ from \refeq{eq:CE-zm}, it is obvious that it
has the required covariance $\vek{C}_{x_a}$ in \refeq{eq:CE-covar-a}.

While \refeq{eq:orth-CE-map}, \refeq{eq:GMKF-RV}, or \refeq{eq:CE-l-filter}
are transformations of $\tnb{x}_f$ to $\tnb{x}_a$ through a simple shift, in
\refeq{eq:CE-l-m-cov} there is an additional linear transformation of the zero mean
part $\vtil{x}^{\vcek{y}}_{a}$.
Although the first two moments of $\vek{x}_{a,\text{cov}}$ in \refeq{eq:CE-l-m-cov}
are correct, it does not seem so simple to proceed further.

In the following, we have two objectives.  For one, we want an assimilated RV
which matches the Bayesian posterior even better, beyond the first two moments,
and on the other hand we do not need so many RVs $\bxi$ and $\bbeta$ to describe
the RV $\vtil{x}^{\vcek{y}}_{a}$.  For the following assume that $\bxi$ and $\bbeta$
are centred  normalised jointly Gaussian and uncorrelated --- hence independent ---
in the \emph{unconditional} expectation $\EXP{\cdot}$.  They are not necessarily
uncorrelated in the conditional expectation $\EXP{\cdot | \vcek{y}}$ though.

To start with this latter task, perform a singular value decomposition on the
$M \times (L+K)$ matrix $\vek{R}$, where normally $L+K \ge M$:
\begin{equation}   \label{eq:zeta-SVD}
  \vek{R} = [\EXP{\vtil{x}^{\vcek{y}}_{a} \otimes \bxi |\vcek{y}},
       \EXP{\vtil{x}^{\vcek{y}}_{a} \otimes \bbeta |\vcek{y}} ]  = 
       \vek{U} \vek{\Sigma} \vek{V}^{\trpos};
\end{equation}
here $\vek{\Sigma}$ is the non-singular diagonal $R\times R$-matrix of non-zero
singular values, with $R\le M$ the \emph{rank} of $\vek{R}$,
and $\vek{U}$ is a $M\times R$ orthogonal matrix ($\vek{U}^{\trpos}\vek{U} = \vek{\Id}_R$),
and $\vek{V}$ a $(L+K)\times R$ orthogonal matrix ($\vek{V}^{\trpos}\vek{V} = \vek{\Id}_R$).
The conditional expectation has to be performed w.r.t.\ $\bxi$ and $\bbeta$
by computing the optimal maps $\vpi_\Psi$ corresponding to $\Psi(\bxi,\bbeta) =
\vtil{x}^{\vcek{y}}_{a} \otimes \bxi$  and $\Psi(\bxi,\bbeta) =
\vtil{x}^{\vcek{y}}_{a} \otimes \bbeta$, denoted by
$\vpi_{\vtil{x} \xi}$ and $\vpi_{\vtil{x} \eta}$ such that
\begin{equation}   \label{eq:tilx-xi-CE}
\EXP{\vtil{x}^{\vcek{y}}_{a} \otimes \bxi |\vcek{y}}=\vpi_{\vtil{x} \xi}(\vcek{y})\qquad\text{and}
\qquad\EXP{\vtil{x}^{\vcek{y}}_{a} \otimes \bbeta |\vcek{y}}=\vpi_{\vtil{x} \eta}(\vcek{y}).
\end{equation}
Similarly, perform an eigenvalue decomposition of the new $\bxi, \bbeta$ correlation matrix
\begin{equation}   \label{eq:xi-eta-corr-EV}
  \vek{C}^{\vcek{y}}_{(\xi \eta)} :=
    \begin{bmatrix}
       \EXP{\bxi \otimes \bxi |\vcek{y}} & \EXP{\bxi \otimes \bbeta |\vcek{y}} \\
       \EXP{\bbeta \otimes \bxi |\vcek{y}}& \EXP{\bbeta \otimes \bbeta |\vcek{y}}
     \end{bmatrix}
   = \begin{bmatrix} \vpi_{\xi \xi}(\vcek{y}) & \vpi_{\xi \eta}(\vcek{y}) \\ 
     \vpi_{\xi \eta}(\vcek{y})^{\trpos} & \vpi_{\eta \eta}(\vcek{y}) \end{bmatrix}
    =  \vek{Q} \vek{\Lambda} \vek{Q}^{\trpos};
\end{equation}
it has size $(L+K)\times(L+K)$, where similarly to \refeq{eq:zeta-SVD}
and \refeq{eq:tilx-xi-CE} above the
conditional expectation has to be performed w.r.t.\ $\bxi$ and $\bbeta$
by computing the optimal maps $\vpi_\Psi$ corresponding to $\Psi(\bxi,\bbeta) =
\bxi \otimes \bbeta$ and the other combinations of $\bxi$ and $\bbeta$, denoted by
$\vpi_{\xi \xi}, \vpi_{\xi \eta}$ and $\vpi_{\eta \eta}$ in \refeq{eq:xi-eta-corr-EV}.
The matrix $\vek{Q}$ of eigenvector-columns is orthogonal
($\vek{Q}^{\trpos}\vek{Q}=\vek{\Id}_{(L+K)}$ and $\vek{Q}\vek{Q}^{\trpos}=\vek{\Id}_{(L+K)}$),
and $\vek{\Lambda}$ is the diagonal matrix of positive eigenvalues.

Then define $R$ new RVs
$\bzeta = [\zeta_1,\dots,\zeta_R]^{\trpos}$ as
\begin{equation}  \label{eq:zeta-def}
   \bzeta =  \vek{V}^{\trpos} \vek{\Lambda}^{-1/2} \vek{Q}^{\trpos}
   \begin{bmatrix} \bxi \\ \bbeta \end{bmatrix} =: \bzeta(\bxi, \bbeta) .
\end{equation}
As a linear transformation of centred jointly Gaussian RVs, the RVs $\bzeta$ are also
centred and  jointly Gaussian, and from \refeq{eq:zeta-def} it is easy to show
that they are also normalised and uncorrelated --- hence independent --- in the
conditional expectation, $\EXP{\bzeta\otimes\bzeta|\vcek{y}} =\vek{\Id}_R$.
In case some eigenvalue vanishes, the matrix $\vek{\Lambda}^{-1/2}$ has to be
understood as the square root of the Moore-Penrose inverse.

A more systematic build-up of the posterior RV beyond the conditional mean
$\vbar{x}^{\vcek{y}}_{a}$ may be achieved with these $R$ new RVs $\bzeta$ as follows.
Choose a set of $J$ linearly independent functions of $\bzeta$:
$\{ \phi_\alpha(\bzeta) \}_{\alpha\in\C{J}_J}$.  As the new RVs $\bzeta$ contain
the same information w.r.t.\ $\vtil{x}^{\vcek{y}}_{a}$ 
as the old set $[\bxi, \bbeta]$, we want to express
$\vtil{x}^{\vcek{y}}_{a}$ in terms of these new RVs.  To this end we form
the Gram matrix $\vek{\Phi}=(\vek{\Phi}_{\alpha \beta})$ with
\begin{equation}   \label{eq:Gram-phi}
\vek{\Phi}_{\alpha \beta} = \EXP{\phi_\alpha(\bzeta(\bxi, \bbeta)) \,
\phi_\beta(\bzeta(\bxi, \bbeta))|\vcek{y}} = \vpi_{\phi_\alpha \phi_\beta}(\vcek{y}),
\end{equation}
where again, as above, the conditional expectation has to be performed
w.r.t.\ $\bxi$ and $\bbeta$ by computing the optimal map
$\vpi_\Psi$ corresponding to $\Psi(\bxi,\bbeta) =
\phi_\alpha(\bzeta(\bxi, \bbeta)) \, \phi_\beta(\bzeta(\bxi, \bbeta))$, denoted by
$\vpi_{\phi_\alpha \phi_\beta}$.

The coefficients in the new expansion
\begin{equation}   \label{eq:FA-xa-zeta}
   \vek{x}^{\vcek{y}}_{a,J}(\bzeta) = \vpi_{x}(\vcek{y}) + 
            \sum_{\alpha\in\C{J}_J} \vek{x}_a^{(\alpha)} \phi_\alpha(\bzeta)
\end{equation}
are then obtained through a Galerkin condition of \refeq{eq:FA-xa-zeta} from
\begin{equation}   \label{eq:xa-coeff}  \forall \alpha\in\C{J}_J:\qquad
   \vek{x}_a^{(\alpha)} = \vek{\Phi}^{-1}\, \vpi_{\vtil{x}\phi_\alpha}(\vcek{y}),
\end{equation}
where again the optimal map $\vpi_{\vtil{x}\phi_\alpha}$ corresponds to the
conditional expectation
\[
    \vpi_{\vtil{x}\phi_\alpha}(\vcek{y}) =
    \EXP{\phi_\alpha(\bzeta(\bxi,\bbeta)) \, \vtil{x}^{\vcek{y}}_{a}(\bxi,\bbeta)|\vcek{y}},
\]
i.e.\ one computes $\vpi_\Psi$ with $\Psi(\bxi,\bbeta) = \phi_\alpha(\bzeta(\bxi,\bbeta)) 
\, \vtil{x}^{\vcek{y}}_{a}(\bxi,\bbeta)$.

As the expression \refeq{eq:FA-xa-zeta} is typically a truncated expansion, the RV
$\vek{x}^{\vcek{y}}_{a,J}(\bzeta)$ will most probably not have the covariance $\vek{C}_{x_a}$
required by \refeq{eq:CE-covar}.  In this case one may use the same procedure as in
\refeq{eq:CE-l-m-cov}.  The covariance of $\vek{x}^{\vcek{y}}_{a,J}$ in \refeq{eq:FA-xa-zeta} is
\begin{equation} \label{eq:CE-covar}  
   \vek{C}_{x_{a,J}} = \EXP{\vtil{x}^{\vcek{y}}_{a,J} \otimes \vtil{x}^{\vcek{y}}_{a,J}|\vcek{y}}
   = \sum_{\alpha, \beta\in\C{J}_J} \vek{\Phi}_{\alpha \beta} 
   \vek{x}_a^{(\alpha)}\otimes\vek{x}_a^{(\beta)} ,
\end{equation}
so that the RV $\vek{x}^{\vcek{y}}_{a,J}(\bzeta)$ in \refeq{eq:FA-xa-zeta} may be corrected
 --- hopefully only slightly --- to
\begin{equation}   \label{eq:FA-xa-cov}
   \vek{x}^{\vcek{y}}_{a,J}(\bzeta) = 
       \vpi_{x}(\vcek{y}) + \vek{C}_{x_a}^{1/2} \vek{C}_{x_{a,J}}^{-1/2}
           \left( \sum_{\alpha\in\C{J}_J} \vek{x}_a^{(\alpha)} \phi_\alpha(\bzeta) \right) .
\end{equation}

\section{Conclusion}         \label{S:conc}
A general approach for state and parameter estimation has been presented in a Bayesian framework.
The Bayesian approach is based here on the conditional expectation (CE) operator,
and different approximations were discussed, where the linear approximation leads to a
generalisation of the well-known Kalman filter (KF), and is here termed the Gauss-Markov-Kalman
filter (GMKF), as it is based on the classical Gauss-Markov theorem.
Based on the CE operator, various approximations
to construct a RV with the proper posterior distribution were shown, where just correcting
for the mean is certainly the simplest type of filter, and also the basis of the GMKF.

Actual numerical computations typically require a discretisation of both the spatial variables
--- something which is practically independent of the considerations here --- and the
stochastic variables.  Classical are sampling methods, but here the use of spectral resp.\
functional approximations is alluded to, and all computations in the examples shown
are carried out with functional approximations.

\bibliography{\thebib/jabbrevlong,\thebib/matthies_BU_paper-1,\thebib/phys_D,\thebib/fa,\thebib/risk,\thebib/stochastics,\thebib/fuq-new,\thebib/sfem,\thebib/highdim,\thebib/filtering,\thebib/inverse}

{ 
   \tiny
       \texttt{\RCSId} 
}

\end{document}